\documentclass[12 pt,reqno]{amsart}
\usepackage{enumitem}
\setlist[enumerate,1]{label={(\alph*)}}

\usepackage{nameref}
\makeatletter
\def\namedlabel#1#2{\begingroup
   #2%
    \def\@currentlabel{#2}%
    \phantomsection\label{#1}\endgroup
}
\makeatother
\usepackage{pdfpages}
\usepackage{amsmath,amsthm,amssymb,amscd}
\usepackage{setspace}
\usepackage{upgreek}

\usepackage{cite}
\usepackage{url}
\usepackage{mathtools}
\usepackage{fullpage}
\usepackage{float}
\usepackage{comment}
\usepackage{longtable}

\usepackage{diagbox}

\usepackage{graphicx}

\usepackage[bookmarksopen,bookmarksdepth=3]{hyperref}

\DeclareFontFamily{U}{rcjhbltx}{}
\DeclareFontShape{U}{rcjhbltx}{m}{n}{<->rcjhbltx}{}
\DeclareSymbolFont{hebrewletters}{U}{rcjhbltx}{m}{n}

\DeclareMathSymbol{\mem}{\mathop}{hebrewletters}{109}%\let\mim\mem
\DeclareMathSymbol{\memh}{\mathop}{hebrewletters}{110}
\DeclareMathSymbol{\nund}{\mathop}{hebrewletters}{141}

\usepackage{etoolbox}
\AtBeginEnvironment{tabular}{\onehalfspacing}

\usepackage{footnotebackref}

\input xy
\xyoption{all}

\usepackage{microtype}

\allowdisplaybreaks[1]

\newtheorem{thm}{Theorem}
\newtheorem{prop}{Proposition}[section]
\newtheorem{lm}[prop]{Lemma}
\newtheorem{cor}[prop]{Corollary}

\theoremstyle{definition}
\newtheorem{dfn}[prop]{Definition}

\theoremstyle{remark}
\newtheorem{rem}[prop]{Remark}

\DeclareMathOperator{\rk}{rk}
\DeclareMathOperator{\Ker}{Ker}
\DeclareMathOperator{\Coker}{Coker}

\DeclareMathOperator{\Id}{Id}

\DeclareMathOperator{\im}{Im}

\DeclareMathOperator{\GW}{GW}
\DeclareMathOperator{\ogw}{OGW}
\DeclareMathOperator{\ogwb}{\overline{\ogw}}

\DeclareMathOperator{\rdim}{rel\,dim}

\newcommand{\rarr}{\rightarrow}
\newcommand{\lrarr}{\longrightarrow}

\newcommand{\R}{\mathbb{R}}
\newcommand{\C}{\mathbb{C}}
\newcommand{\Z}{\mathbb{Z}}

\newcommand{\M}{\mathcal{M}}

\renewcommand{\P}{\mathbb{CP}}
\renewcommand{\L}{\Lambda}

\newcommand{\g}{\Gamma}

\newcommand{\Hh}{\widehat{H}}

\newcommand{\p}{\mathfrak{p}}

\newcommand{\q}{\mathfrak{q}}

\renewcommand{\d}{\partial}

\renewcommand{\Im}{\im}

\newcommand{\Ob}{\overline{\Omega}}

\newcommand{\RP}{\mathbb{RP}}

\newcommand{\oFS}{\omega_{FS}}

\newcommand{\sly}{\Pi}

\renewcommand{\a}{\alpha}

\newcommand{\Ah}{\widehat{A}}

\renewcommand{\k}{\Bbbk}
\newcommand{\Lc}{\L_c}

\newcommand{\s}{\mathfrak{s}}
\newcommand{\hp}{\mathcal{H}}
\newcommand{\ml}{\mathcal{L}}
\newcommand{\mla}{\ml^\star}
\newcommand{\md}{{m^*}}
\newcommand{\la}{l^\star}

\renewcommand{\k}{\Bbbk}
\renewcommand{\i}{i}
\newcommand{\bP}{\mathbb{P}}

\makeatletter
\DeclareRobustCommand{\cev}[1]{%
  \mathpalette\do@cev{#1}%
}
\newcommand{\do@cev}[2]{%
  \fix@cev{#1}{+}%
  \reflectbox{$\m@th#1\vec{\reflectbox{$\fix@cev{#1}{-}\m@th#1#2\fix@cev{#1}{+}$}}$}%
  \fix@cev{#1}{-}%
}
\newcommand{\fix@cev}[2]{%
  \ifx#1\displaystyle
    \mkern#23mu
  \else
    \ifx#1\textstyle
      \mkern#23mu
    \else
      \ifx#1\scriptstyle
        \mkern#22mu
      \else
        \mkern#22mu
      \fi
    \fi
  \fi
}

\makeatother
\usepackage{xpatch}

\makeatletter
\xpretocmd{\@adminfootnotes}{\let\@makefntext\BHFN@OldMakefntext}{}{}
\renewcommand\@makefntext[1]{%
  \@ifundefined{@makefnmark}
    {}
    {%
     \renewcommand\@makefnmark{%
       \mbox{%
         \textsuperscript{%
           \normalfont
           \hyperref[\BackrefFootnoteTag]{\@thefnmark}%
         }%
       }\,%
     }%
     \BHFN@OldMakefntext{#1}%
  }%
}
\makeatother

\DeclareMathOperator{\spa}{Span}

\newcommand{\altfrac}[2]{\ifmmode\def\tmp{$}\else\def\tmp{}\fi\mbox{%
    {\raisebox{.24\ht\strutbox}{\tmp#1\tmp}}%
    \kern-2.2pt\scalebox{1.6}[1.5]{/}\kern-1.8pt%
    {\tmp#2\tmp}%
    }}

\newcommand{\dl}{\Delta}
\renewcommand{\cir}{L}

\newcommand{\os}{\mathfrak{a}}
\newcommand{\At}{{\widetilde{A}}}
\newcommand{\Ht}{{\widetilde{H}}}
\renewcommand{\cir}{L}

\newcommand{\ups}{\Upsilon}
\newcommand{\pwr}{\delta}

\newcommand{\uu}{\mathbf{u}}

\newcommand{\Lh}{\widehat{L}}

\newcommand{\E}{\mathcal{E}}

\newcommand{\Xh}{\widehat{X}}

\usepackage{marginnote}

\title{Examples of relative quantum cohomology}

\author[K. Hugtenburg]{Kai Hugtenburg}
\address{School of Mathematics and Statistics\\ Lancaster University \\Bailrigg\\Lancaster\\ LA1 4YW \\UK}
\email{kai.hugtenburg@gmail.com}
\author[S. Tukachinsky]{Sara B. Tukachinsky}
\address{School of Mathematical Sciences\\ Tel Aviv University\\Tel Aviv, 6997801, Israel }\email{sarabt1@gmail.com}

\date{Feb. 2024}
\subjclass[2020]{53D45, 14N35 (Primary) 53D37, 53D12 (Secondary)}
\begin{document}

\begin{abstract}
We compute the quantum cohomology relative to a Lagrangian submanifold in some complete intersections. For quadric hypersurfaces, we also give a full computation of the genus zero open Gromov-Witten invariants.
\end{abstract}

\maketitle
\setcounter{tocdepth}{1}
\tableofcontents

\section{Introduction}
Let $(X,\omega)$ be a symplectic manifold of dimension $2n,$ and let $L \subset X$ be a relatively spin Lagrangian submanifold.
In~\cite{ST3}, the relative quantum cohomology $QH(X,L)$ is defined under certain conditions, for example, if $L$ is a real cohomology sphere. The associativity of the relative quantum product is shown to be equivalent to the open and closed WDVV equations.

This paper gives explicit computations of the small relative quantum ring structure for some complete intersections, as well as an explicit treatment and a full computation of all open Gromov-Witten invariants of the quadric, as follows.

We start by formulating the results for Beauville-type~\cite{Beauville} complete intersections.
Let
$X\subset \P^{n+r}$ with $n\ge 3$ be a complete intersection of degree $(d_1,\ldots, d_r)$.
The first Chern class of $X$ is given by
\[
c_1(X)= n+1-\sum_{j=1}^r(d_j-1).
\]
For the purposes of this section, assume that $n > 2\sum_{j=1}^r(d_j-1)-1$, meaning,
\begin{equation}\label{eq:c1bd}
2c_1(X)> n+1.
\end{equation}
Throughout, we use the notation
\[
\ups:=\prod_{j=1}^rd_j^{d_j},\qquad \pwr:=\sum_{j=1}^r(d_j-1),
\qquad d:=\prod_{j=1}^rd_j.
\]
Let $v_1,\ldots,v_p\in H^n(X;\R)$ be the primitive classes. For all $j,k\in\{1,\ldots,p\}$, write
\[
g_{jk}:=\int_X v_j\wedge v_k.
\]
Finally, we will consider the following polynomials in the formal variables $a,b,c$:
\begin{gather*}
r(a):=a^{n+1} - q\ups\cdot a^{\pwr}%\label{rel:Rx}
,\\
r_{jk}(a,b,c):=bc - \frac{1}{d}g_{jk}\cdot \big(a^n-q \ups \cdot a^{\pwr-1}\big)%\label{rel:Rjk}
.
\end{gather*}
The following is immediate from the main theorem of~\cite{Beauville}.
Here $u$ is a variable that corresponds to the hyperplane class in $H^2(X;\R)$, and $v_1,\ldots,v_p,$ correspond to the primitive classes in $H^n(X;\R)$.
Then the absolute quantum cohomology is given by
\begin{equation}\label{eq:beauville}
QH^*(X) \simeq
\altfrac{\R[[q]][u,v_1,\ldots,v_p]}{\big(r(u), r_{jk}(u,v_j,v_k), uv_j)_{j,k\in\{1,\ldots,p\}}}.
\end{equation}

Let $L\subset X$ be a Lagrangian submanifold such that $L$ is homeomorphic to $S^n$.
In particular, $\varpi =\Id$ and we take the unique relative spin structure with $w_\s=0$.
In view of Lemma~\ref{lm:mvp} below, we generalize the assumption~\eqref{eq:c1bd} to
\begin{equation}\label{eq:mubd}
\mu > n+1.
\end{equation}

Taking $x$ to be a variable that corresponds to the hyperplane class in $\Hh^2(X,L;\R)$ and let $w_1,\ldots,w_P,$ correspond to the primitive classes in $\Hh^n(X,L;\R)$. In particular, $\{w_j\}_{j=1}^P$ map to $\{v_j\}_{j=1}^p$ under $\Hh^*(X,L;\R)\to H^*(X)$.

Consider the case $[L]\ne 0\in H_n(X;\R)$. Then $P=p-1$, and we may assume without loss of generality that $v_p$ corresponds to the Poincar\'e dual of $[L]$. Then we prove the following.

\begin{thm}\label{thm:ciLneq0}
If $[L]\ne 0$, then the relative quantum cohomology is given by
\[
QH^*(X,L)\simeq \altfrac{\R[[q]][x,w_1,\ldots,w_{p-1}]}{\big(r(x), r_{jk}(x,w_j,w_k), xw_j\big)_{j,k\in\{1,\ldots,p-1\}}}.
\]
The relative quantum cohomology injects into the absolute one via
\begin{gather*}
\os :QH^*(X,L)\longrightarrow QH^*(X)\\
q\mapsto q,\quad x\mapsto u, \quad w_j\mapsto v_j.
\end{gather*}
\end{thm}

\begin{rem}
One can verify directly by inspecting the relations ideal that the formula for $QH^*(X,L)$ indeed defines a subalgebra. Similarly, in the next theorem we get a a rank-1 extension of the algebra $QH^*(X)$. This is a general phenomenon, explained in more detail in Sections~\ref{ssec:rings} and~\ref{sssec:rqh} below.
\end{rem}

Now consider the case $[L] = 0\in H_n(X;\R)$. Then $P=p$, and $\Hh^*(X,L;\R)$ has an additional generator to which we associate the variable $y$.
Then we prove the following.

\begin{thm}\label{thm:ciL0}
If $[L]=0$, then the relative quantum cohomology is given by
\[
QH^*(X,L)\simeq
\altfrac{\R[[q]][x,w_1,\ldots,w_p,y]}{\big(r(x), r_{jk}(x,w_j,w_k), xw_j, y^2, yx, yw_j\big)_{j,k\in\{1,\ldots,p\}}}.
\]
The relative quantum cohomology surjects onto the absolute one via
\begin{gather*}
\os :QH^*(X,L)\longrightarrow QH^*(X)\\
q\mapsto q,\quad x\mapsto u,\quad w_j \mapsto v_j,\quad y\mapsto 0.
\end{gather*}
\end{thm}

For our first example computed explicitly,
consider
$(X,L) = (\C P^n,\R P^n).$ It is shown in~\cite{RuanTian, KontsevichManin} that the absolute quantum cohomology of $\C P^n$ is given by the ring $\altfrac{\R[[q]][u]}{(u^{n+1}-q)}.$ To allow the absolute quantum cohomology to interact with the relative quantum cohomology, we extend the Novikov ring by adjoining $q^{1/2}$, and we twist the signs of quantum cohomology by the class in $H^2(\C P^n;\Z/2\Z)$ determined by the relative spin structure on $\R P^n.$ Thus we take
\[
QH^*(\C P^n)\simeq \altfrac{\R[[q^{1/2}]][u]}{(u^{n+1}-(-1)^{\frac{n+1}{2}}q)}.
\]
\begin{thm}\label{thm:rqhcpn}
The relative quantum cohomology of the odd-dimensional projective space is given by
\[
QH^*(\C P^n,\R P^n)\simeq
\altfrac{\R[[q^{1/2}]][x,y]}{I}
\]
with
\[
I=
(y^2-2q^{1/2}y,\quad
x^{n+1}-(-1)^{\frac{n+1}{2}}q
-(-1)^{\frac{n-1}{2}}\frac{1}{2}q^{1/2}y, \quad xy = 0).
\]
The relative quantum cohomology surjects onto the absolute one via
\begin{gather*}
\os :QH^*(\P^n,\RP^n)\longrightarrow QH^*(\P^n)\\
q^{1/2}\mapsto q^{1/2}, \quad x\mapsto u, \quad y\mapsto 0.
\end{gather*}
\end{thm}

\begin{rem}
Note that the case $(\P^n,\RP^n)$ is not covered by the Beauville-style complete intersection theorem, since $\mu=n+1$. This indeed makes the relations in the quantum ring substantially different: While in Theorem~\ref{thm:ciL0} the additional generator $y$ gave rise to an independent nilpotent component, here $y$ satisfies non-trivial relations, and in particular changes non-trivially the relation for $x^{n+1}$.
\end{rem}

Further, we carry out direct, definition-based computations of the small relative quantum cohomology in the special case when $X$ is a quadric hypersurfaces and $L \simeq S^n$ is the real locus with respect to a particular real structure for arbitrary $n > 1$. Again, we take the unique relative spin structure with $w_\s=0$. This family of examples, while in some ways straightforward, offers a comprehensive overview of the principles at play. Indeed, the Lagrangian submanifold of the quadric is homologically nontrivial if and only if the quadric is even dimensional.
Direct calculations in the two cases, of being homologically trivial or nontrivial, require the different tools of the wall-crossing formula and the open-closed relation, respectively, proven in~\cite{ST3} and cited below as Theorems~\ref{thm:wallcross} and~\ref{thm:magic}. More concretely, we get the following results.

Let $X$ be a quadric hypersurface in $\P^{n+1}$ with odd $n$ given by $\sum_{j=0}^{n}z_j^2 - z_{n+1}^2=0$. It is shown in~\cite{FultonPandharipande} for $n=3$ and in~\cite{Beauville} for $n\ge 3$ that the absolute quantum cohomology of $X$ is given by
\[
QH^*(X)\simeq \altfrac{\R[[q]][u]}{(u^{n+1}-4qu)}.
\]
See also Proposition~\ref{prop:qhquadod} below.
\begin{thm}\label{thm:rqhqodd}
Let $L\subset X$ be the Lagrangian submanifold given by the real locus of complex conjugation.
The relative quantum cohomology is given by
\[
QH^*(X,L)\simeq
\altfrac{\R[[q]][x,y]}{I}
\]
with
\[
I=(x^{n+1}-4qx,\quad xy,\quad y^2).
\]
The relative quantum cohomology surjects onto the absolute one via
\begin{gather*}
\os :QH^*(X,L)\longrightarrow QH^*(X)\\
q\mapsto q, \quad x\mapsto u, \quad y\mapsto 0.
\end{gather*}
\end{thm}

Now, let $X$ be a quadric hypersurface in $\P^{n+1}$ with even $n$ given by $\sum_{j=0}^{n}z_j^2 - z_{n+1}^2=0$. As shown in~\cite{Beauville} for $n\ge 4$ and revisited in Proposition~\ref{prop:qhquadev} below,

the absolute quantum cohomology of $X$ is given by
\[
QH^*(X)\simeq \altfrac{\R[[q]][u,v]}{I},
\]
with
\[
I=(uv, \quad v^2- (-1)^{n/2+1}4q-(-1)^{n/2}u^n).
\]
\begin{thm}\label{thm:rqhqeven}
Let $L\subset X$ be the Lagrangian submanifold given by the real locus of complex conjugation.
The relative quantum cohomology is given by
\[
QH^*(X,L)\simeq \altfrac{\R[[q]][x]}{(x^{n+1}-4qx)}.
\]
The relative quantum cohomology injects into the absolute one via
\begin{gather*}
\os :QH^*(X,L)\longrightarrow QH^*(X)\\
q\mapsto q,\quad x\mapsto u.
\end{gather*}
\end{thm}

Lastly, let $X$ be the quadric hypersurface in $\P^3$ given by $\sum_{j=0}^2z_j^2-z_3^2=0$. As shown in~\cite[Example 7.2]{CrauderMiranda},
the absolute quantum cohomology of $X$ is given by
$
\altfrac{\R[[q_1,q_2]][l,\la]}{(l^2-q_2, {\la}^2-q_1)}.
$
As with the example of the projective space, we need to modify the Novikov ring in order to allow the absolute quantum cohomology to interact with the relative one. As explained in Section~\ref{sec:quad2}, this is done by identifying $q_1$ and $q_2$. Thus we take
\[
QH^*(X)\simeq \altfrac{\R[[q]][l,\la]}{(l^2-q, \quad {\la}^2-q)}.
\]
Changing basis to $u = l+\la$ and $v = l - \la$ gives
\[
QH^*(X)\simeq \altfrac{\R[[q]][u,v]}{(u^2+v^2-4q, \quad uv)},
\]
Which agrees with the result obtained for the general even-dimensional quadric hypersurfaces.

\begin{thm}\label{thm:rqhq2}
Let $L\subset X$ be the Lagrangian submanifold given by the real locus of complex conjugation.
The relative quantum cohomology is given by
\[
QH^*(X,L)\simeq \altfrac{\R[[q]][x]}{(x^3-4qx)}.
\]
The relative quantum cohomology injects into the absolute one via
\begin{gather*}
\os :QH^*(X,L)\longrightarrow QH^*(X)\\
q\mapsto q,\quad x\mapsto u.
\end{gather*}
\end{thm}

For each of the above explicit examples we verify that all the closed and open Gromov-Witten invariants are computable via WDVV-based recursions. This is akin to spelling out the big quantum product. The computability of the closed invariants of the odd dimensional quadric is given by~\cite{KontsevichManin, RuanTian}, while the computation of the open invariants is new.
For the even dimensional quadric, even the closed invariants require special care, since the cohomology of the ambient manifold is not generated by divisors.
We prove the following:

\begin{thm}
	\label{thm:oddcomputability}
Let $X$ be a quadric hypersurface in $\P^{n+1}$ with odd $n$ given by $\sum_{j=0}^{n}z_j^2 - z_{n+1}^2=0$. Let $L \cong S^n$ be the real locus of $X$. The open Gromov-Witten invariants of $(X,L)$ are entirely determined by the open WDVV equations, the axioms of $\ogw$, the wall-crossing formula Theorem \ref{thm:wallcross}, the closed Gromov-Witten invariants of $X$, and $\ogw_{1,3} = (-1)^{\frac{(n-2)(n-3)}{2}}2$.
\end{thm}
We observe an interesting property of the open Gromov-Witten invariants of the odd dimensional quadrics, which follows from a general result regarding compatibility of the open WDVV equations and sign changes (see Proposition \ref{prop: sign change OWDVV}).
\begin{prop}
	For the odd dimensional quadrics, suppose that we change the sign of the initial value $\ogw_{1,3}$, then all the other $\ogw$ only change by a sign computed in Proposition \ref{prop: sign change odd quadrics}.
\end{prop}

We include some sample values of open Gromov-Witten invariants for various odd quadrics (tables \ref{tb:ogwoddboundaryonly} and \ref{tb:ogw323}).
\begin{table}[ht]
	\centering
\begin{tabular}{l|lll}
	\diagbox{$n$}{$m$} & $0$ & $1$ & $2$ \\\hline
	$3$ & $2$ & $96$ & $1036800$ \\
	$5$ & $-2$ & $11520$ & $-2807046144000$ \\
	$7$ & $2$ & $2580480$ & $26876909061144576000$ \\
	$9$ & $-2$ & $928972800$ & $-684658186021829527732224000$ \\
	$11$ & $2$ & $490497638400$ & $38607476855046755233605322014720000$ \\
\end{tabular}
	\caption{Values of $OGW^n_{\beta,k}$ with $\beta = 1 + (n-1)m$ and $k = 3 + 2nm$.}\label{tb:ogwoddboundaryonly}
\end{table}

\begin{table}[ht]
	\centering
	\begin{tabular}{l|lllll}
		\diagbox{$l_3$}{$\beta$} & $1$ & $2$ & $3$ & $4$ & $5$ \\\hline
		$0$ & $2$ & $24$ & $11744$ & $16767744$ & $65095308288$ \\
		$1$ & $2$ & $4$ & $1488$ & $1501248$ & $4557798144$ \\
		$2$ & $0$ & $0$ & $208$ & $143296$ & $336418560$ \\
		$3$ & $0$ & $0$ & $32$ & $14624$ & $26285312$ \\
		$4$ & $0$ & $0$ & $0$ & $1600$ & $2185728$ \\
		$5$ & $0$ & $0$ & $0$ & $192$ & $194560$ \\
	\end{tabular}
	\caption{Values of $OGW^3_{\beta,1}(\Gamma_2^{\otimes l_2}\otimes\Gamma_3^{\otimes l_3})$ with $l_2 = 3\beta-1-2l_3$.}\label{tb:ogw323}
\end{table}

\begin{thm}
	\label{thm:evencomputability}
Let $X$ be a quadric hypersurface in $\P^{n+1}$ with even $n$ given by $\sum_{j=0}^{n}z_j^2 - z_{n+1}^2=0$. Let $L \cong S^n$ be the real locus of $X$. The open Gromov-Witten invariants of $(X,L)$ are entirely determined by the axioms of $OGW$, the closed Gromov-Witten invariants of $X$, and the open-closed relation Theorem \ref{thm:magic}. Moreover, the closed Gromov-Witten invariants of $X$ are entirely determined by the axioms of $GW$, the WDVV equations, and $GW_1(\Delta_{n-1},\Delta_n) = 4$.
\end{thm}
\begin{thm}
	\label{thm:q2computability}
	Let $X$ be a quadric hypersurface in $\P^{3}$ given by $\sum_{j=0}^{3} z_j^2 - z_4^2 = 0$. Let $L \cong S^2$ be the real locus of $X$. The open Gromov-Witten invariants of $(X,L)$ are entirely determined by the axioms of $OGW$, the closed Gromov-Witten invariants of $X$, and the open-closed relation Theorem \ref{thm:magic}. Moreover, the closed Gromov-Witten invariants of $X$ are entirely determined by the axioms of $GW$, the WDVV equations, and the basic closed Gromov-Witten invariants (see Lemma \ref{lm:gw2nn-1}).
\end{thm}
We include some sample values of open Gromov-Witten invariants for various odd quadrics (tables \ref{tb:ogw423} and \ref{tb:ogw42L}).

\begin{table}[ht]
	\centering
	\begin{tabular}{l|llll}
		\diagbox{$l_3$}{$\beta$} & $1$ & $2$ & $3$ & $4$ \\\hline
		$0$ & $2$ & $412$ & $3186280$ & $137299159040$ \\
		$1$ & $2$ & $116$ & $564024$ & $17481735264$ \\
		$2$ & $0$ & $32$ & $102704$ & $2279623168$ \\
		$3$ & $0$ & $8$ & $18992$ & $303098752$ \\
		$4$ & $0$ & $0$ & $3456$ & $40795008$ \\
	\end{tabular}
	\caption{Values of $OGW^4_{\beta,1}(\Gamma_2^{\otimes l_2}\otimes\Gamma_3^{\otimes l_3})$ with $l_2 = 4\beta-1-2l_3$.}\label{tb:ogw423}
\end{table}

\begin{table}[ht]
	\centering
	\begin{tabular}{l|lll}
		\diagbox{$m$}{$\beta$} & $1$ & $2$ & $3$ \\\hline
		$0$ & $-2$ & $44$ & $-14760$ \\
		$1$ & $2$ & $-4$ & $5544$ \\
		$2$ & $0$ & $-52$ & $-9160$ \\
		$3$ & $0$ & $412$ & $42120$ \\
		$4$ & $0$ & $0$ & $-284904$ \\
	\end{tabular}
	\caption{Values of $OGW^4_{\beta,1}(\Gamma_2^{\otimes l_2}\otimes PD(L)^{\otimes l_L})$ with $l_2 = 1+ 2m$ and $l_L = 4\beta - 2 - 2m$.}\label{tb:ogw42L}
\end{table}

\subsection{Acknowledgments}
The authors would like to thank
Jake Solomon for essentially collaborating on parts of this manuscript, particularly Section~\ref{sec:quadtop} and the appendix.
The first author would like to thank Nick Sheridan for suggesting the project of computing open Gromov-Witten invariants for odd-dimensional quadrics.
The first author was partly supported by the Royal Society Research Grant for Research Fellows:
RGFnR1n181009ERC, the ERC Starting Grant 850713 – HMS and EPSRC Grant EP/W015749/1. The second author was partly supported by NSF grant No. DMS-1638352, ISF grant No. 2793/21, and the Colton Foundation.

\section{Background}

\subsection{Setting and notation}

Let $(X,\omega)$ be closed symplectic manifold of dimension $2n$, $J$ an $\omega$-tame almost complex structure, and $L\subset X$ a Lagrangian submanifold with a relative spin structure $\s$.
Recall that $\s$ determines a class
\[
w_{\s} \in H^2(X;\Z/2\Z)
\]
such that $w_2(TL) = i^* w_{\s}$. By abuse of notation we think of $w_{\s}$ as a homomorphism $w_\s:H_2(X;\Z)\to \Z/2\Z$.
Let
\[
\varpi:H_2(X;\Z)\lrarr H_2(X,L;\Z)
\]
be the map from the long exact sequence of the pair $(X,L)$.
Let
\[
\mu: H_2(X,L;\Z) \lrarr \Z
\]
be the Maslov index as defined in~\cite{CieliebakGoldstein}. The following is well known.
\begin{lm}\label{lm:mvp}
We have $\mu(\varpi(\beta))= 2c_1(\beta).$
\end{lm}

\subsection{Bases}\label{ssec:rings}

Consider the subcomplex of differential forms on $X$ consisting of those with trivial integral on $L$,
\[
\Ah^*(X,L):=\left\{\eta\in A^*(X)\;\bigg|\,\int_L\eta=0\right\},
\]
and write $\Hh^*(X,L;\R):=H^*(\Ah^*(X,L), d)$.
This version of relative cohomology fits in the long exact sequence
\begin{equation}\label{eq:yrho}
\xymatrix{
\Hh^*(X,L;\R) \ar[rr]^\rho && H^*(X;\R) \ar[dl]^{\i_\R} \\
 & \R[-n] \ar[ul]_{[1]}^y,
}
\end{equation}
where
\[
\i_\R: H^*(X;\R) \lrarr \R
\]
is the map given by integration over $L$.
If $[L]\ne 0\in H_n(X;\R)$,
the long exact sequence~\eqref{eq:yrho} becomes a short exact sequence
\begin{equation}\label{eq:splitLn0} 0\lrarr \Hh^*(X,L)\stackrel{\rho}{\lrarr}H^*(X)\stackrel{i_\R}{\lrarr}\R[-n] \lrarr 0.
\end{equation}
If $[L]=0$, the sequence becomes
\begin{equation}\label{eq:splitL0}
0\lrarr \R[-n-1]\stackrel{y[-1]}{\lrarr} \Hh^*(X,L;\R) \stackrel{\rho}{\lrarr} H^*(X;\R)\lrarr 0.
\end{equation}
Set $M:=\rk H^*(X;\R)-1$, $N:=\rk \Hh^*(X,L;\R)-1$, and $K= \min \{M,N\}$. In particular, $K +1 = \max\{M,N\}.$
Choose homogeneous bases
\[
\dl_0,\ldots,\dl_M\in H^*(X;\R)\quad \text{and} \quad\g_0,\ldots,\g_N\in \Hh^*(X,L;\R)
\]
such that $\rho(\g_j)=\dl_j$ for $j=0,\ldots, K,$ and if $N > M$ also $\rho(\g_N) = 0$. In particular, $\deg \dl_j=\deg \g_j$ for $j\le K$.
Additionally, set
\begin{equation}\label{eq:dL}
\g_\diamond:=y(1)\in \Hh^*(X,L;\R)\quad \text{and} \quad \dl_L:=PD([L])\in H^*(X;\R).
\end{equation}
Note that exactly one of the two does not vanish, depending on whether $[L]\in H_*(X;\R)$ is trivial or not.
If $[L]=0$, we assume without loss of generality that the bases are chosen so that $\g_N=\g_\diamond$.
For convenience, we define $\g_M:=0$ if $M>N$ and $\dl_N:=0$ if $N>M$.

\begin{rem}
By~\cite[Lemma~5.14]{ST2}, $\Hh^*(X,L;\R)$ is isomorphic to the standard (reduced) relative cohomology when $H^*(L;\R)\simeq H^*(S^n;\R)$. In particular, this is true in all the cases considered in the current paper.
\end{rem}

\subsection{Rings}

Define a Novikov coefficient ring
\begin{gather*}
\L=\left\{\sum_{i=0}^\infty a_iT^{\beta_i}\bigg|a_i\in\R,\beta_i\in H_2(X,L;\Z),\omega(\beta_i)\ge 0,\; \lim_{i\to \infty}\omega(\beta_i)=\infty\right\}.
\end{gather*}
Grading on $\L$ is defined by declaring $T^\beta$ to be of degree $\mu(\beta).$

Let $\g_0,\ldots,\g_N,$ and $\dl_0,\ldots,\dl_M,$ be the bases from Section~\ref{ssec:rings}.
Let $t_0,\ldots,t_{K+1},$ be formal variables with degrees $\deg t_j=2-\dl_j$ for $j\le M$ and $\deg t_j = 2-\deg \g_j$ for $j\le N$. Let $s$ be a formal variable with degree $\deg s = 1-n$. We will use $\dl_j,\g_j,$ for interior constraints of closed and open Gromov-Witten invariants, respectively, while the variable $s$ will correspond to boundary point constraints for open invariants.

Define graded-commutative rings
\[
R:=\L[[t_0,\ldots,t_N,s]],\qquad Q:=\L[[t_0,\ldots,t_N]],
\qquad Q_X:=\L[[t_0,\ldots,t_M]].
\]
Define
\[
\rho^* :Q_X\lrarr Q
\]
to be the $\L$ algebra homomorphism determined by
\[
\rho^*(t_j) = \begin{cases}
t_j, & j\le K,\\
0, & j>K.
\end{cases}
\]
The ring $R$ comes with a valuation
\[
\nu:R\lrarr \R,
\]
given by
\[
\nu\left(\sum_{j=0}^\infty a_jT^{\beta_j}s^k\prod_{i=0}^Nt_i^{l_{ij}}\right)
= \inf_{\substack{j\\a_j\ne 0}} \left(\omega(\beta_j)+k+\sum_{i=0}^N l_{ij}\right),
\]
and similarly for $Q$ and $Q_X$.
For $\Upsilon\in \{R, Q, Q_X\}$, write $H^*(X;\Upsilon)$ and $\Hh^*(X,L;\Upsilon)$ for $H^*(X;\R)\otimes_{\R} \Upsilon$ and $\Hh^*(X,L;\R)\otimes_{\R} \Upsilon$, respectively. Here the tensor is completed with respect to the valuation on $\Upsilon$.

\begin{rem}
Comparing with the notation convention of~\cite{ST3}, throughout the paper we take $\sly=H_2(X,L;\Z)$ and $U=H^*(X)$. Furthermore, here we write $R$, $Q$, and $Q_X$ instead of $R_W$, $Q_W$, and $Q_U$, respectively; we write $QH^*(X)$ for $QH^*_U(X)$ and $QH^*(X,L)$ for $QH_U(X,L)$.
\end{rem}
\begin{rem}
In~\cite{ST3}, we used two different Novikov coefficient rings, $\L$ and $\Lc$.
This was essential in the construction of a cone complex that provides a natural setting for definitions and proofs in open Gromov-Witten theory, where all cases can be treated uniformly.
The current paper, however, has more limited aims, which allows us to simplify exposition by removing $\Lc$.
\end{rem}

\subsection{Closed GW}
The genus zero closed Gromov-Witten invariants are a collection of maps
\[
GW_\beta : H^*(X;\R)^{\otimes l}\longrightarrow \R, \qquad\quad \beta \in H_2(X;\Z), \quad l \in \Z_{\geq 0}.
\]
Intuitively, for $A_1,\ldots,A_l\in H^*(X;\R)$, the invariant $\GW_\beta(A_1,\ldots,A_l)$ is the algebraic count of $J$-holomorphic spheres of degree $\beta$ that pass through chains representing the Poincar\'e duals of the classes $A_j$.

These invariants can be combined into the potential function $\Phi\in Q_X$, which can be written explicitly in coordinates as follows.
Let $\dl_j\in H^*(X;\R)$ be the basis from Section~\ref{ssec:rings} associated with the variables $t_0,\ldots,t_M$. Then
\[
\Phi = \Phi(\{t_i\}_{i =0}^M)
=
\sum_{\substack{\beta\in H_2(X;\Z)\\ r_i\ge 0}}
(-1)^{w_\s(\beta)}\frac{T^{\varpi(\beta)}\prod_{i=0}^M t_{M-i}^{r_{M-i}}}{\prod_{i=0}^M r_i!} \GW_{\beta}(\otimes_{i =0}^M\dl_i^{\otimes r_i}).
\]
The non-standard sign $(-1)^{w_\s(\beta)}$ allows the closed Gromov-Witten invariants encoded in $\Phi$ to interact with the open Gromov-Witten invariants of the Lagrangian $L,$ which depend on the relative spin structure $\s.$ See in particular Theorem~\ref{thm:OWDVV} below.

As shown in, e.g., \cite{KontsevichManin, MS, RuanTian}, the Gromov-Witten invariants have the following properties.
\begin{thm}[Closed Gromov-Witten axioms]
\label{axioms_c}
The invariants $\GW$ of $X$ satisfy the following axioms.
Let $A_j \in H^*(X;\R)$ for $j = 1,\ldots,l$.
\begin{enumerate}
	\item (Degree)
		$\GW_{\beta}(A_1,\ldots,A_l)=0$ unless
		\begin{equation}\label{ax_c_deg}
		2n-6+2c_1(\beta)+2l = \sum_{j=1}^l |A_j|.
		\end{equation}
	\item (Unit / Fundamental class)
		\begin{equation}\label{ax_c_unit}
		\GW_{\beta}(1,A_{1},\ldots,A_{l-1})=
		\begin{cases}
			\int_XA_1\wedge A_2, & (\beta,l)=(\beta_0,3),\\
			0, & \text{otherwise}.
		\end{cases}
		\end{equation}
	\item (Zero)
		\begin{equation}\label{ax_c_zero}
		\GW_{\beta_0}(A_1,\ldots,A_l)=
		\begin{cases}
		\int_X A_1\wedge A_2\wedge A_3, & l=3,\\
		0, & \text{otherwise}.
		\end{cases}
		\end{equation}
	\item (Divisor)
		If $|A_l| = 2,$ then
		\begin{equation}\label{ax_c_divisor}
		\GW_{\beta}(A_1,\ldots,A_{l})=\int_\beta A_l \cdot\GW_{\beta}(A_1,\ldots,A_{l-1}).
		\end{equation}
\end{enumerate}
\end{thm}
In addition, the potential $\Phi$ satisfies the Witten-Dijkgraaf-Verlinde-Verlinde (WDVV) equations, given in \cite[Theorem-Definition 4.5]{KontsevichManin} and~\cite[Theorem 3]{RuanTian0}. To formulate them, we need some additional notation.
Set
\[
g_{ij}:=\int_X \dl_i\wedge \dl_j,
\]
and define $(g^{ij})_{i,j=0}^M$ to be the matrix inverse to $(g_{ij})_{i,j=0}^M$.
\begin{thm}[Closed WDVV equations]\label{thm:WDVV}
Let $u,v,w\in \{0,\ldots,M\}$.
Then,
\begin{equation}\label{eq:wdvv}
\sum_{m,l=0}^M\d_{t_u}\d_{t_v}\d_{t_m}\Phi \cdot g^{ml}\cdot \d_{t_l}\d_{t_w}\d_{t_y}\Phi=
(-1)^{|\Delta_u|(|\Delta_v| + |\Delta_w|)}\sum_{m,l=0}^M\d_{t_v}\d_{t_{w}}\d_{t_m}\Phi\cdot g^{ml} \cdot\d_{t_l}\d_{t_{u}}\d_{t_y}\Phi.
\end{equation}
\end{thm}

\subsection{Open GW}
In many situations, genus zero open Gromov-Witten invariants are defined in~\cite{ST2,ST3}.
For the purposes of this paper, we focus on the situation when $H^*(L;\R)\simeq H^*(S^n;\R)$, which includes all the examples discussed here.
The genus zero open Gromov-Witten invariants are a collection of maps
\[
\ogw_{\beta,k} : \Hh^*(X,L;\R)^{\otimes l} \to \R, \qquad \beta \in H_2(X,L;\Z), \quad k,l \in \Z_{\geq 0}.
\]
Intuitively, for $A_1,\ldots,A_l\in \Hh^*(X,L;\R)$, the invariant $\ogw_{\beta,k}(A_1,\ldots, A_l)$ counts configurations of $J$-holomorphic disks that collectively have degree $\beta$ and pass through $k$ points in $L$ and chains in $X\setminus L$ representing the Poincar\'e duals of $A_j$; which disks join into such configurations is determined by a bounding chain. See~\cite[Section 1.2.4]{ST2} for more details. The open Gromov-Witten invariants can be arranged in the open superpotential $\Omega\in R$, which can be written explicitly in coordinates as follows.
Let $\g_j\in \Hh^*(X,L;\R)$ be the basis from Section~\ref{ssec:rings} associated with the variables $t_0,\ldots,t_N,$ and think of the additional variable $s$ as corresponding to the Poincar\'e dual of a point in $H^n(L;\R)$.
Then the superpotential $\Omega\in R$ is given by
\begin{equation}\label{eq:O}
\Omega = \Omega(s,\{t_i\}_{i =0}^N)=\sum_{\substack{\beta\in\sly\\k\ge 0\\r_i\ge 0}}\frac{T^\beta s^k \prod_{i =0}^Nt_{N-i}^{r_{N-i}}}{k! \prod_{i=0}^Nr_i!} \ogw_{\beta,k}(\otimes_{i =0}^N\g_i^{\otimes r_i}).
\end{equation}
By design, if $k = 0$ and $\beta\in \im \varpi$ then $\ogw_{\beta,k}\equiv 0$. See~\cite[Remark 1.8]{ST2}. In the case when $[L]=0\in H_n(X;\R)$, enhanced invariants
\[
\ogwb_{\beta,k} : \Hh^*(X,L;\R)^{\otimes l} \to \R, \qquad \beta \in H_2(X,L;\Z), \quad k,l \in \Z_{\geq 0},
\]
are defined in~\cite{ST3}, such that
\begin{equation}\label{eq:ogw}
\ogw_{\beta,k}(\eta_1,\ldots,\eta_l) =
\begin{cases}
0, & k = 0 \text{ and } \beta \in \im\varpi, \\
\ogwb_{\beta,k}(\eta_1,\ldots,\eta_l), & \text{otherwise}.
\end{cases}
\end{equation}
These can be arranged into the enhanced superpotential
\begin{equation}\label{eq:Ob}
\Ob(s,\{t_i\}_{i = 0}^N)=\sum_{\substack{\beta\in\sly\\k\ge 0\\r_i\ge 0}}\frac{T^\beta s^k \prod_{i=0}^Nt_{N-i}^{r_{N-i}}}{k! \prod_{i=0}^Nr_i!} \ogwb_{\beta,k}(\otimes_{i=0}^N\g_i^{\otimes r_i}).
\end{equation}
When $k = 0$ and $\beta\in \im \varpi$, the enhanced invariants $\ogwb_{\beta,k}$ do not necessarily vanish. However, their values depend on an additional piece of data: A homomorphism
\begin{equation}\label{eq:PR}
P_\R: \Hh^*(X,L;\R)\lrarr \Coker(\i_\R)
\end{equation}
that is a left inverse to the map $\bar y : \Coker \i_\R \to \Hh^*(X,L;\R)$ induced by $y$ from the long exact sequence~\eqref{eq:yrho}.
When $[L]\ne 0$, note that $\Coker(\i_\R)=0$ and therefore $P_\R\equiv 0$.
When $[L]=0$, a choice of $P_\R$ is equivalent to a choice of a splitting in~\eqref{eq:splitL0}. Another perspective on the choice of $P_\R$ is explained in~\cite[Remark 4.12]{ST3}: It is equivalent to the choice of a homology class $V\in H_{n+1}(X,L;\R)$ such that $\d V=[L]$.

The open Gromov-Witten invariants satisfy the following properties analogous to the closed ones.
For $\ogw$ this is~\cite[Theorem 4]{ST2}, and for $\ogwb$ this is~\cite[Proposition 4.19]{ST3}.
\begin{thm}[Open Gromov-Witten axioms]
\label{axioms}
The invariants $\ogwb$ of $(X,L)$ satisfy the following axioms.
Let $A_j \in W$ for $j = 1,\ldots,l$.
\begin{enumerate}
	\item (Degree)
		$\ogwb_{\beta,k}(A_1,\ldots,A_l)=0$ unless
		\begin{equation}\label{ax_deg}
		n-3+\mu(\beta)+k+2l = kn+ \sum_{j=1}^l |A_j|.
		\end{equation}
	\item (Unit / Fundamental class)
		\begin{equation}\label{ax_unit}
		\ogwb_{\beta,k}(1,A_{1},\ldots,A_{l-1})=
		\begin{cases}
			-1, & (\beta,k,l)=(\beta_0,1,1),\\
            P_\R(A_1), & (\beta,k,l) = (\beta_0,0,2),\\
			0, & \text{otherwise}.
		\end{cases}
		\end{equation}
	\item (Zero)
		\begin{equation}\label{ax_zero}
		\ogwb_{\beta_0,k}(A_1,\ldots,A_l)=
		\begin{cases}
		-1, & (k,l)=(1,1)\text{ and } A_1=1,\\
        P_\R(A_1 \smallsmile A_2), & (k,l) = (0,2),\\
		0, & \text{otherwise}.
		\end{cases}
		\end{equation}
	\item (Divisor)
		If $|A_l| = 2,$ then
		\begin{equation}\label{ax_divisor}
		\ogwb_{\beta,k}(A_1,\ldots,A_{l})=\int_\beta A_l \cdot\ogwb_{\beta,k}(A_1,\ldots,A_{l-1}).
		\end{equation}
\end{enumerate}
The invariants $\ogw$ satisfy the same axioms with $P_\R$ replaced with $0$.
\end{thm}

The superpotential $\Omega$ satisfies an open version of the WDVV equations.
To formulate it, we need some additional notation.
Recall that $M=\rk H^*(X;\R)-1,$ $N=\rk \Hh^*(X,L;\R)-1$, and $K=\min\{M,N\}$. In particular, $K=\rk \Ker(P_\R)-1$.
Possibly modifying our choice of $\{\g_j\}_{j = 0}^K$ from Section~\ref{ssec:rings}, we can assume that
$\{\g_j\}_{j=0}^K\subset \Ker(P_{\R})$
is a basis.
Set
\[
g_{ij}:=\int_X \dl_i\wedge \dl_j,
\qquad i,j=0,\ldots, M,
\]
and define $(g^{ij})_{i,j=0}^{M}$ to be the matrix inverse to $(g_{ij})_{i,j=0}^{M}$.
Below is the statement of~\cite[Corollary 1.6]{ST3}.

\begin{thm}[Open WDVV equations for point-like bounding chains]\label{cor:owdvv}\label{thm:OWDVV}
Let $u,v,w\in \{0,\ldots,N\}$.
If $[L]=0,$ then
\begin{multline}\label{eq:cor1}
\sum_{\substack{0\le l \le K\\ 0\le m\le M}}\d_{t_u}\d_{t_l}\Ob\cdot g^{lm}\cdot
\rho^*\d_{t_m}\d_{t_w}\d_{t_v}\Phi
-\d_{t_u}\d_s\Ob\cdot\d_{t_w}\d_{t_v}\Ob=\\
=
\sum_{\substack{0\le l \le M\\ 0\le m\le K}}\rho^*\d_{t_u}\d_{t_w}\d_{t_l}\Phi \cdot g^{lm}\cdot
\d_{t_m}\d_{t_v}\Ob
-\d_{t_u}\d_{t_w}\Ob\cdot\d_{t_v}\d_s\Ob,
\end{multline}
and
\begin{equation}\label{eq:cor2}
\sum_{\substack{0\le l \le K\\ 0\le m\le M}}\d_s\d_{t_l}\Ob\cdot g^{lm}\cdot
\rho^*\d_{t_m}\d_{t_w}\d_{t_v}\Phi
-\d_s^2\Ob\cdot\d_{t_w}\d_{t_v}\Ob=
-\d_s\d_{t_w}\Ob\cdot\d_{t_v}\d_s\Ob.
\end{equation}
If $[L]\ne 0,$ then equation~\eqref{eq:cor1} holds with $\Ob$ replaced by $\Omega.$
\end{thm}

\subsection{Interplay}
In addition to the standard Gromov-Witten-type axioms, we have two phenomena unique to the open setup.
The first is a wall-crossing property, which allows us to trade boundary constraints for interior ones if $L$ is homologically trivial in $X$.
This is~\cite[Theorem 6]{ST3}:

\begin{thm}[Wall-crossing]\label{thm:wallcross}
Suppose $[L] = 0$.
The invariants $\ogwb$ satisfy
\begin{equation}\label{eq:wallcross}
\ogwb_{\beta,k+1}(\eta_1,\ldots,\eta_l)
=
-\ogwb_{\beta,k}(\Gamma_\diamond,\eta_1,\ldots,\eta_l).
\end{equation}
\end{thm}

The second relates open invariants with closed ones when $L$ is homologically nontrivial in~$X$.
This is~\cite[Corollary 1.5]{ST3}:

\begin{thm}[Open-closed relation]\label{thm:magic}
Suppose that $[L] \neq 0$.
Let $\eta_1,\ldots,\eta_l \in \Hh^*(X,L).$ If $k\ge 2$, or if $k\ge 1$ and $\beta \notin \Im \varpi$, then
\[
\ogwb_{\beta,k}(\eta_1,\ldots,\eta_l)=0.
\]
Moreover, for $\beta\in \Im(\varpi)$ and $\eta \in H^*(X;\R)$ such that $\int_L \eta = 1,$
\begin{equation}\label{eq:magic}
\ogwb_{\beta,1} (\eta_1,\ldots\eta_l)
=
\sum_{\substack{\hat\beta \in H_2(X;\Z)\\  \varpi(\hat\beta)=\beta}}
(-1)^{n+1+w_\s(\hat\beta)}
\GW_{\hat\beta}(\eta,\dl_L,\rho(\eta_1), \ldots,\rho(\eta_l)).
\end{equation}
\end{thm}

\subsection{Sign changes of the open GW}
\label{section: sign changes}
We also observe the following result about compatibility of the open Gromov-Witten axioms (Theorem \ref{axioms}) and the open WDVV equations (Theorem \ref{cor:owdvv}) with specific sign changes. When all open Gromov-Witten invariants are determined from certain initial values by the axioms and the open WDVV equations, this sometimes allows us to deduce how the invariants change when the initial values change sign. See Proposition \ref{prop: sign change odd quadrics}. For the cleanest statements, we need to consider invariants with interior constraints only from $\ker(P_\R)$
\begin{dfn}
	Restrict the invariants $\ogwb$ (or $\ogw$ when $[L] \neq 0$) to invariants \begin{equation*}
		\ogwb'_{\beta,k}: \ker(P_\R)^{\otimes l} \rightarrow \R
	\end{equation*}
	and consider the restricted superpotential \begin{equation}
		\Ob'(s,\{t_i\}_{i = 0}^{N-1})=\sum_{\substack{\beta\in\sly\\k\ge 0\\r_i\ge 0}}\frac{T^\beta s^k \prod_{i=0}^{N-1}t_{N-i}^{r_{N-i}}}{k! \prod_{i=0}^{N-1}r_i!} \ogwb_{\beta,k}(\otimes_{i=0}^{N-1}\g_i^{\otimes r_i}).
	\end{equation}
	Equivalently, $\Ob' := \Ob|_{t_N = 0}$.
\end{dfn}
The following is clear from the definition:
\begin{lm}
	The invariants $\ogwb'$ and $\ogw'$ satisfy the open Gromov-Witten axioms with $P_\R$ replaced with $0$.
\end{lm}
We also have: \begin{lm}
	The restricted superpotentials $\Ob'$ and $\Omega'$ still satisfy the open WDVV equations.
\end{lm}
\begin{proof}
	In the open WDVV equations, the summation indices run over $0 \leq l \leq K$ and $0 \leq m \leq M$ (as any integral $\int_X \Gamma_\diamond \Gamma_j = 0$ for any $j$). This ensures that the equations still hold after restriction to $t_N = 0$.
\end{proof}
\begin{dfn}
	\label{dfn: sign change}
	For an arbitrary integer $r$ consider the sign-changed invariants \begin{equation*}
		\ogwb^\circ_{\beta,k} := (-1)^{\frac{k-1}{r}}\ogwb'_{\beta,k},
	\end{equation*}
	with a similar definition with $\ogw$ in place of $\ogw$. We then also form the sign-changed superpotentials $\Ob^\circ$ and $\Omega^\circ$ by using the invariants $\ogwb^\circ$ and $\ogw^\circ$ in their definitions.
\end{dfn}
\begin{rem}
	Of course, in order for the invariants to be real numbers (or even rational) the above only makes sense if $\ogwb'_{\beta,k} = 0$ for all $k$ such that $r \nmid k-1$. We will ignore this issue for now.
\end{rem}
Finally, we can state the new result:
\begin{prop}
	\label{prop: sign change OWDVV}
	For any choice of $r$, the invariants $\ogwb^\circ$ and $\ogw^\circ$ still satisfy the open Gromov-Witten axioms, and the superpotentials $\Ob^\circ$ and $\Omega^\circ$ still satisfy the open WDVV equations.
\end{prop}
\begin{proof}
	The first statement holds as the only non-zero invariants featuring in the open Gromov-Witten axioms have $k = 1$. For the second statement, consider the coefficient of $s^k$ in Equation \eqref{eq:cor1}. The first term on the left-hand side has its sign changed by $(-1)^{\frac{k-1}{r}}$ compared to the equation for $\Ob$. The second term has it changed by $(-1)^{\frac{k_1}{r}}(-1)^{\frac{k_2-1}{r}} = (-1)^{\frac{k-1}{r}}$ as required. Similar reasoning holds for the right-side of the equation, so that the overall equality is maintained. The same holds for Equation \eqref{eq:cor2}.
\end{proof}

\subsection{Closed quantum cohomology}

Write $QH^*(X):=H^*(X;\L)$.
The small quantum product
\[
*: QH^*(X)\otimes QH^*(X)\lrarr QH^*(X)
\]
was first introduced in the physics literature~\cite{LercheVafaWarner, Witten1}, and given a mathematical definition in~\cite{RuanTian0,RuanTian, KontsevichManin}.
The version we use in the current work is given by the formula
\begin{equation}\label{eq:dqp}
\dl_v*\dl_u =  \sum_{\substack{0\le l,m\le M\\ d\in H_2(X;\Z)}} (-1)^{w_\s(d)}T^{\varpi(d)}
\GW_d(\dl_v,\dl_u,\dl_l)\cdot g^{lm}\cdot\dl_m, \qquad 0 \leq u,v \leq M.
\end{equation}
This product is well known~\cite{RuanTian0,RuanTian, KontsevichManin} to be associative and graded commutative.

\subsection{Relative quantum cohomology}\label{sssec:rqh}

The small relative quantum product is defined analogously, though the definition is a bit more delicate. Write $QH^*(X,L):=\Hh^*(X,L;\L)$.
Then the product
\footnote{Read $\mem$ as ``mem''.}
\[
\mem: QH^*(X,L)\otimes QH^*(X,L) \lrarr QH^*(X,L)
\]
is given by
\begin{align}
\mem(\g_v, \g_u) 
= &
\sum_{\substack{0\le m\le K\\0\le l\le M\\ d\in H_2(X;\Z)}}
(-1)^{w_\s(d)}T^{\varpi(d)}
\GW_d(\dl_v,\dl_u,\dl_l) \cdot g^{lm}\cdot\g_m+ \label{eq:decomp}\\
&\qquad+\sum_{d\in H_2(X,L;\Z)}T^d\ogwb_{d,0}(\g_v,\g_u)\cdot \g_\diamond, & 0 \leq u,v \leq N. \notag
\end{align}
In the special case $[L]\ne 0$, since $\g_\diamond=0$, the expression simplifies to
\begin{align*}
\mem(\g_v, \g_u)
= &
\sum_{\substack{0\le m\le K\\ 0\le l\le M\\d\in H_2(X;\Z)}}
(-1)^{w_\s(d)}T^{\varpi(d)}
\GW_d(\dl_v,\dl_u,\dl_l) \cdot g^{lm}\cdot\g_m, \qquad 0 \leq u,v \leq N.
\end{align*}
The product $\mem$ is shown in~\cite{ST3} to be associative and graded commutative.
Though in~\cite{ST3} the definition of $\mem$ is based on a chain-level construction, the above explicit formulas follow from~\cite[Lemma 5.5]{ST3}.
\begin{rem}
The definition of $\mem$ given in~\cite{ST3} is free of choices. The explicit expression~\eqref{eq:decomp} a priori appears to depend on the choice of $P_{\R}$ through the definition of $\ogwb$, but this cancels with the dependence on $P_{\R}$ through the condition that $\{\g_j\}_{j=0}^K\subset \Ker(P_{\R})$ is a basis.
\end{rem}

\begin{lm}\label{lm:homom}
The map $\rho$ from~\eqref{eq:yrho} induces an algebra homomorphism
\[
\os: QH^*(X,L)\lrarr QH^*(X).
\]
\end{lm}
\begin{proof}
This follows from equation~(53) of~\cite{ST3}.
\end{proof}

\subsection{Big closed quantum cohomology}

The product $*$ can be deformed by $\dl=\sum_{j=0}^Mt_j\dl_j$, resulting in the big closed quantum product. More explicitly,
set
$QH_{big}^*(X):=QH^*(X)\otimes Q$.
Then
\[
*^{big}:QH_{big}^*(X)\otimes QH_{big}^*(X)\lrarr QH_{big}^*(X)
\]
is given by
\[
\dl_v *^{big}\dl_u =
\sum_{\substack{0\le l,m\le M\\ d\in H_2(X;\Z)\\ p\ge 0 }} (-1)^{w_\s(d)}T^{\varpi(d)}\frac{1}{p!}\GW_d(\dl^{\otimes p}\otimes \dl_v \otimes \dl_u \otimes\dl_l)\cdot g^{lm}\cdot\dl_m,
\qquad 0\le u,v,\le M.
\]

This product is well known~\cite{RuanTian0,RuanTian, KontsevichManin} to be associative and graded commutative.

\subsection{Big relative quantum cohomology}

The product $\mem$ can be deformed by $\g=\sum_{j=0}^Nt_j\g_j$, resulting in the big relative quantum product. More explicitly, set $QH^*_{big}(X,L):=QH^*(X,L)\otimes Q$. The product
\[
{\mem}^{big}:QH^*_{big}(X,L)\otimes QH^*_{big}(X,L)\lrarr QH^*_{big}(X,L)
\]is given by
\begin{align*}
{\mem}^{big}(\g_v ,\g_u)
&=
\sum_{\substack{0\le m\le K\\0\le l\le M\\d\in H_2(X;\Z)\\ p\ge 0}}
(-1)^{w_\s(d)}T^{\varpi(d)}
\frac{1}{p!}\GW_d(\dl^{\otimes p}\otimes \dl_v\otimes\dl_u\otimes\dl_l) \cdot g^{lm}\cdot\g_m+\\
&\qquad+\sum_{\substack{d\in H_2(X,L;\Z)\\p\ge 0}}T^d\frac{1}{p!}\ogwb_{d,0}(\g^{\otimes p}\otimes\g_v\otimes\g_u)\cdot \g_\diamond,
\hspace{5em}0 \leq u,v \leq N.
\end{align*}
When $[L]\ne 0$ the formula simplifies to
\begin{multline*}
{\mem}^{big}(\g_v, \g_u)
=
\sum_{\substack{0\le m\le K\\ 0\le l\le M\\d\in H_2(X;\Z)\\ p\ge 0}}
(-1)^{w_\s(d)}T^{\varpi(d)}
\frac{1}{p!}\GW_d(\dl^{\otimes p}\otimes \dl_v\otimes\dl_u\otimes\dl_l) \cdot g^{lm}\cdot\g_m,\\
 0 \leq u,v \leq N.
\end{multline*}

%As with the big closed quantum cohomology, writing the relative ring in terms of generators and relations is generally not feasible. However, the product of any two elements can be computed directly from the definition if all the values of $\GW$ and $\ogwb$ are known.

The product $\mem^{big}$ is shown in~\cite{ST3} to be associative and graded commutative.
Though in~\cite{ST3} the definition of $\mem^{big}$ is based on a chain-level construction, the above explicit formulas follow from~\cite[Lemma 5.5]{ST3}.
\begin{rem}
The definition of $\mem^{big}$ given in~\cite{ST3} is free of choices. The explicit expression~\eqref{eq:decomp} a priori appears to depend on the choice of $P_{\R}$ through the definition of $\ogwb$, but this cancels with the dependence on $P_{\R}$ through the condition that $\{\g_j\}_{j=0}^K\subset \Ker(P_{\R})$ is a basis.
\end{rem}

\begin{lm}
The map $\rho$ from~\eqref{eq:yrho} induces an algebra homomorphism
\[
\os_{big}: QH^*_{big}(X,L)\lrarr QH^*_{big}(X).
\]
\end{lm}
\begin{proof}
This follows from equation~(53) of~\cite{ST3}.
\end{proof}

\section{Complete intersections}

In this section, we prove Theorems~\ref{thm:ciLneq0} and~\ref{thm:ciL0}.

\subsection{Closed quantum cohomology}

We choose a basis $\{\Delta_j\}_{j=0}^M$ of $H^*(X;\R)$ as follows.
Let $\Delta_1\in H^2(X;\R)$ be the hyperplane class,
namely, writing $\oFS$ for the Fubini-Study form on $\P^{n+r}$ and taking
\begin{equation}\label{eq:incl}
\iota:X\hookrightarrow \P^{n+r}
\end{equation}
to be the inclusion, we take
$\Delta_1=[\iota^*\oFS]\in H^2(X;\R)$.
Let $\Delta_j:=\Delta_1^{\wedge j}$ for $j\in\{0,1,\ldots,n\}$. Let $\{\Delta_j\}_{j=n+1}^{M}$ be generators of the primitive cohomology such if $[L]\ne  0$ then $\Delta_M=PD([L])$.
Further, assume without loss of generality that $\Delta_j$ for $j>n$ are chosen so that $\{\Delta_{j}\}_{j=0}^K$ span $\Ker(i_{\R})$. Denote the number of primitive basis elements by
\[
p:= M-n.
\]
To choose a basis for the relative cohomology $\Hh^*(X,L;\R)$, and subsequently computing the quantum cohomology, we have to discuss the cases of $[L]$ being trivial and non-trivial separately.

\subsection{\texorpdfstring{$L$}{L} is homologically non-trivial}

We choose a basis $\{\g_j\}_{j=0}^N$ of $\Hh^*(X,L;\R)$ as follows.
Let $\g_1\in H^2(X,L;\R)$ be the hyperplane class, namely, for $\iota$ the inclusion from~\eqref{eq:incl}, we set $\g_1=[\iota^*\oFS]\in \Hh^2(X,L;\R)$.
For $j\in\{0,1,\ldots,n\}$, let $\g_j:=\g_1^{\wedge j}$.
In particular, $\rho(\g_j)=\Delta_j$.
Let $\{\g_j\}_{j=n+1}^{N}$ be generators of the primitive part of $\Hh^*(X,L;\R)$, namely, homogeneous classes independent of $\g_{n/2}$ (when $n$ is even) such that $\rho(\g_j)=\Delta_j$ for $j=n+1,\ldots,N$. Note that when $[L]\ne 0$ we have $N=M-1$, and recall that $\Delta_M = PD([L])$.

\begin{proof}[Proof of Theorem~\ref{thm:ciLneq0}]
Lemma~\ref{lm:homom}, exact sequence~\eqref{eq:yrho} and the assumption that $L$ is homologically non-trivial, imply that $QH^*(X,L)$ is the subalgebra of $QH(X,L)$ given by $\ker(\i_\R)\otimes \R[[q]].$
\end{proof}

\subsection{\texorpdfstring{$L$}{L} is homologically trivial}

We choose a basis $\{\g_j\}_{j=0}^N$ of $\Hh^*(X,L;\R)$ as follows.
Let $\g_1\in H^2(X,L;\R)$ be the hyperplane class, namely, for $\iota$ the inclusion from~\eqref{eq:incl}, we set $\g_1=[\iota^*\oFS]\in \Hh^2(X,L;\R)$.
For $j\in\{0,1,\ldots,n\}$, let $\g_j:=\g_1^{\wedge j}$.
In particular, $\rho(\g_j)=\Delta_j$.
Let $\{\g_j\}_{j=n+1}^{N-1}$ be generators of the primitive part of $\Hh^*(X,L;\R)$, namely, homogeneous classes independent of $\g_{n/2}$ (when $n$ is even) such that $\rho(\g_j)=\Delta_j$ for $j=n+1,\ldots,N-1$. Note that when $[L]= 0$ we have $N-1=M$. Finally, set $\g_N :=\g_\diamond$ from~\eqref{eq:dL}.

\begin{lm}\label{lm:degvan}
For any homogeneous $a_1,a_2\in \Hh^*(X,L)$ with $|a_1|+|a_2|\le 2n+2$, we have $\ogwb_{\beta,0}(a_1,a_2)=0$.
\end{lm}
\begin{proof}
By the degree axiom~\eqref{ax_deg},
the value of $\ogwb_{\beta,0}(a_1,a_2)$ can only be nonzero if
\[
n+1+\mu(\beta)=|a_1|+|a_2| \le 2n+2.
\]
But then $\mu(\beta)\le n+1$, contradicting the assumption~\eqref{eq:mubd} and Lemma~\ref{lm:mvp}.
\end{proof}

\begin{proof}[Proof of Theorem~\ref{thm:ciL0}]
Let $a\in \{y,x,w_1,\ldots,w_p\}$.
Since $y$ is trivial in absolute cohomology, the $\GW$ term in equation~\eqref{eq:decomp} vanishes in $\mem(y,a)$ and we obtain
\[
\mem(y,a)=\sum_{\beta\in H_2(X,L;\Z)} T^\beta \ogwb_{\beta,0}(y,a)\cdot y.
\]
By Lemma~\ref{lm:degvan}, the right-hand side vanishes. Thus,
\begin{equation}\label{eq:ya}
ya =0 \quad \mbox{ for }\quad a\in \{y,x,w_1,\ldots,w_p\}.
\end{equation}

Recall the short exact sequence~\eqref{eq:splitL0}. Let
\[
\sigma: H^*(X;\R)\lrarr \Hh^*(X,L;\R)
\]
be the section of $\rho$ induced by $P_{\R}$ from~\eqref{eq:PR}. In particular, $\sigma(u)=x$ and $\sigma(v_j)=w_j$.
For
\[
(a_1,a_2)\in \{(x,x^j)\,|\,1\le j\le n\} \cup \{(x,w_j), (w_j,w_k) \,|\, j,k\in \{1,\ldots,p\}\},
\]
equation~\eqref{eq:decomp} and Lemma~\ref{lm:degvan} give
\begin{align}
\mem(a_1, a_2)
= &
\sum_{\substack{0\le l,m\le K=M\\ d\in \Z_{\ge 0}}}
T^{d}
\GW_d(\rho(a_1),\rho(a_2),\dl_l) \cdot g^{lm}\cdot\g_m
+\sum_{d\in \Z_{\ge 0}}T^d\ogwb_{d,0}(a_1,a_2)\cdot \g_\diamond \notag\\
= &
\sum_{\substack{0\le l,m\le K=M\\ d\in \Z_{\ge 0}}}
T^{d}
\GW_d(\rho(a_1),\rho(a_2),\dl_l) \cdot g^{lm}\cdot\g_m \notag\\
=& \sigma(\rho(a_1)*\rho(a_2)). \label{eq:a1a2}
\end{align}
By induction, we obtain
\begin{equation}\label{eq:x^i}
x^i = \sigma(u^i),\qquad i = 1,\ldots,n.
\end{equation}
Combining equations~\eqref{eq:a1a2},~\eqref{eq:x^i}, and~\eqref{eq:beauville}, we obtain
\begin{equation}\label{eq:xn1}
x^{n+1} = \mem(x,x^n) = \sigma(u*u^n) = \sigma(q\Upsilon u^\delta) = q\Upsilon \sigma(u^\delta) = q\Upsilon x^\delta.
\end{equation}
Similarly,
\begin{equation}\label{eq:wjkpf}
w_jw_k =\sigma(v_jv_k)
= \frac{1}{d}g_{jk}\cdot \big(\sigma(u^n-q \ups \cdot u^{\pwr-1})\big)
= \frac{1}{d}g_{jk}\cdot \big(x^n-q \ups \cdot x^{\pwr-1}\big).
\end{equation}
Lastly,
it follows directly from~\eqref{eq:a1a2} and~\eqref{eq:beauville} that
\[
xw_j=0,\qquad j=1,\ldots,p.
\]
Lemma~\ref{lm:homom} and equation~\eqref{eq:beauville} imply that $1,x,\ldots,x^n,w_1,\ldots,w_p,y,$ generate $QH(X,L)$ as a $\L$-module. From the associativity of $\mem$ it follows
that~\eqref{eq:ya},~\eqref{eq:xn1}, and~\eqref{eq:wjkpf},
are the only non-trivial relations in $QH(X,L)$.
\end{proof}

\section{Projective space}\label{sec:cpn}
\subsection{Topological preliminaries}
In this section, $(X,L)=(\P^n,\RP^n)$ for odd $n\geq 3.$ We take $\omega = \oFS$ the Fubini-Study form.
We equip $\RP^n\subset \P^n$ with the relative spin structure that gives
\[
\ogw_{1,2}=2,
\]
as explained in~\cite[Lemma 6.1]{ST3}.
We write $\g_j = [\omega^j] \in \Hh^*(X,L;\R)$ and $\dl_j = [\omega^j] \in H^*(X;\R)$, for $j=0,\ldots,n$. Then $\{\dl_j\}_{j=0}^n$ form a basis of $H^*(X;\R)$. Since $[L]=0\in H_n(X;\Z)$, by short exact sequence~\eqref{eq:splitL0}, the classes $\{\g_j\}_{j=0}^n\cup \{\g_\diamond=y(1)\}$ form a basis of $\Hh^*(X,L;\R)$.

\begin{lm}\label{lm:rcp}
\begin{enumerate}
\item
$H_2(X,L;\Z) \simeq \Z.$
\item \label{it:vpcp}
The map $\varpi : H_2(X;\Z) \to H_2(X,L;\Z)$ is given by multiplication by $2.$
\end{enumerate}
\end{lm}

Identify $H_2(X,L;\Z)$ with $\Z$ so that the non-negative integers correspond to classes $\beta \in H_2(X,L;\Z)$ with $\omega(\beta) \geq 0$, and similarly for $H_2(X;\Z)$.
It follows from Lemma~\ref{lm:rcp} that for any $\beta\in H_2(X,L;\Z)$ and $\hat\beta\in H_2(X;\Z)$, we have
\begin{equation}\label{eq:divpn}
\int_\beta \g_1=\frac{1}{2}\beta, \qquad
\int_{\hat\beta} \dl_1 = \hat\beta.
\end{equation}
These values will be used when applying the divisor axioms~\eqref{ax_divisor} and~\eqref{ax_c_divisor}.

\begin{lm}\label{lm:charp}
The relevant characteristic classes of $(X,L)$ are as follows:
\begin{enumerate}
\item\label{it:c1p}
For all $\hat\beta\in H_2(X;\Z)$ we have
\[
c_1(\hat\beta)=(n+1)\hat\beta,\qquad
w_\s(\hat\beta)\equiv \frac{n+1}{2}\hat\beta \pmod 2.
\]
\item\label{it:mup}
For all $\beta\in H_2(X,L;\Z)$, we have
$
\mu(\beta)=(n+1)\beta.
$
\end{enumerate}
\end{lm}

For $i,j\in\{0,\ldots,n\}$, the pairing on absolute cohomology is given by
\begin{equation}\label{eq:gijp}
g_{ij}= \int_X \dl_i\wedge\dl_j =\delta_{i,n-j},
\end{equation}
where $\delta_{k,l}$ is the Kronecker delta.
Therefore the inverse matrix $(g^{ij})_{i,j}:={(g_{ij})_{i,j}}^{-1}$ is also given by
$g^{ij}=\delta_{i,n-j}$.

\subsection{Absolute quantum ring}

Using equation~\eqref{eq:gijp}, Lemmas~\ref{lm:rcp} and~\ref{lm:charp}\ref{it:c1p}, and the definition of the quantum product~\eqref{eq:dqp}, for $i,j\in\{0,\ldots,n\}$, we obtain
\[
\dl_i*\dl_j=\sum_{\substack{0\le m\le n\\ d\ge 0}} (-1)^{\frac{n+1}{2}\cdot d}T^{2d}\GW_d(\dl_i,\dl_j,\dl_m)\dl_{n-m}.
\]

It is shown in~\cite[Section 5.3]{KontsevichManin} and~\cite[Example 8.5]{RuanTian} that, over the coefficient ring $\R[[q]]$, the absolute quantum cohomology of the projective space is given by $\altfrac{\R[[q]][x]}{(x^{n+1}-q)}$. Here, $x$ corresponds to $\dl_1$ and $q$ corresponds to $T^{\varpi(1)}=T^2$. For our purposes, we modify the setup as follows. First, we extend the coefficient ring by adjoining $q^{1/2}$, which corresponds to $T=T^1$, thereby obtaining the coefficient ring $\R[[q^{1/2}]]\simeq\L$. Further, we twist $\GW_{\hat\beta}$ by the sign $(-1)^{w_\s(\hat\beta)}=(-1)^{\frac{n+1}{2}\hat\beta}$. Thus, we have the following.

\begin{prop}
The absolute quantum cohomology is given by
\[
QH^*(X)\simeq \altfrac{\R[[q^{1/2}]][x]}{( x^{n+1}-(-1)^{\frac{n+1}{2}}q)}.
\]
\end{prop}

\subsection{Relative quantum ring}
Using equation~\eqref{eq:gijp}, Lemmas~\ref{lm:rcp} and~\ref{lm:charp}\ref{it:c1p}, and the definition of the relative quantum product~\eqref{eq:decomp},
for $i,j\in \{0,\ldots,n\}\cup\{\diamond\}$, we obtain
\begin{equation}\label{eq:smempn}
\mem(\g_i,\g_j)
= \sum_{\substack{0\le m\le n\\ d\ge 0}} (-1)^{\frac{n+1}{2}\cdot d}T^{2d}\GW_d(\dl_i,\dl_j,\dl_m)\g_{n-m}
+ \sum_{d\ge 0}T^d\ogwb_{d,0}(\g_i,\g_j)\g_\diamond.
\end{equation}

The following is included in~\cite[Theorem 10]{ST3}.
\begin{lm}\label{lm:pninit}
$\ogw_{1,2}=2$, $\ogw_{1,1}(\g_{\frac{n+1}{2}})=0$, $\ogw_{1,0}(\g_n)=(-1)^{\frac{n-1}{2}}$.
\end{lm}

\begin{proof}[Proof of Theorem~\ref{thm:rqhcpn}]
By the degree axiom for $\GW$~\eqref{ax_c_deg} and $\ogwb$~\eqref{ax_deg}, all contributions with $d\ge 2$ to the sums in~\eqref{eq:smempn} are zero.
Consider now invariants of the form $\ogwb_{d,0}(\g_i,\g_j)$, constituting the coefficients in the second sum in~\eqref{eq:smempn}. For $d=0$, the zero axiom~\eqref{ax_zero} implies that the only nonzero values are when $i=\diamond$ and $j=0$ or vice versa. For $d=1$, the degree axiom~\eqref{ax_deg} implies that nonzero values are only possible when $|\g_i|+|\g_j|=2n+2$.

For $i,j\in\{0,\ldots,n\}$, this means there is potentially contribution if $j=n+1-i$, that is, an invariant of the type $\ogwb_{1,0}(\g_i,\g_{n+1-i})$.
By Lemma~\ref{lm:pninit}, the wall crossing theorem~\eqref{eq:wallcross}, zero axiom~\eqref{ax_zero}, and divisor axiom~\eqref{ax_divisor}, the other possibly nontrivial values of $\ogwb_{d,0}(\g_i,\g_j)$ for $d=0,1,$ are explicitly computed to be
\begin{gather*}
\ogwb_{0,0}(\g_{\diamond},1)=1, \\
\ogwb_{1,0}(\g_{\diamond},\g_{\diamond}) =\ogwb_{1,2}=2,\\
\ogwb_{1,0}(\g_{\diamond},\g_{\frac{n+1}{2}}) =-\ogwb_{1,1}(\g_{\frac{n+1}{2}})=0,\\
\ogwb_{1,0}(\g_1,\g_n) =\frac{1}{2}\ogwb_{1,0}(\g_n)=\frac{1}{2}\cdot(-1)^{\frac{n-1}{2}}.
\end{gather*}
This immediately implies that
\[
\mem(\g_1,\g_\diamond)=0.
\]
Further, combining the above computation of $\ogwb_{d,0}(\g_i,\g_j)$ with the standard computation of the closed small quantum product, we deduce the following relations:
\begin{align*}
&\mem(\g_{\diamond},1)=\g_{\diamond},\\
&\mem(\g_{\diamond},\g_{\diamond})= 2T\g_\diamond,\\
&\mem(\g_\diamond,\g_{\frac{n+1}{2}}) = 0,\\
&\mem(\g_1,\g_n)= (-1)^{\frac{n+1}{2}}T^2+(-1)^{\frac{n-1}{2}}\frac{1}{2}T\g_\diamond,\\
&\mem(\g_i,\g_j)=\g_{i+j},\qquad i,j\in\{0,\ldots,n\},\quad i+j\le n.
\end{align*}
From the associativity of $\mem$ we see that $\ogwb_{1,0}(\g_i,\g_j) = \ogwb_{1,0}(\g_1,\g_n)$, which does not generate new relations.
Taking $q,x,y,$ to be the variables corresponding to $T^2,\g_1,\g_\diamond,$ respectively, we get the desired result.
\end{proof}

\begin{rem}
As noted above, it follows from the associativity of $\mem$ that
\[
\mem(\g_{\frac{n+1}{2}},\g_{\frac{n+1}{2}})
= \mem(\g_1,\g_n)
= (-1)^{\frac{n+1}{2}}T^2+(-1)^{\frac{n-1}{2}}\frac{1}{2}T\g_\diamond.
\]
This implies that
\[
\ogwb_{1,0}(\g_{\frac{n+1}{2}},\g_{\frac{n+1}{2}})=(-1)^{\frac{n-1}{2}}\frac{1}{2}.
\]
This corrects a sign error made in Section 1.3.11 of~\cite{ST3}.
\end{rem}

\subsection{Full computation of invariants}

\begin{prop}
All $\GW$ invariants are determined by the WDVV equations and the initial value $\GW_{1}(\dl_n^{\otimes 2})=1$.
All $\ogw$ invariants are determined by the OWDVV equations, the wall-crossing theorem~\eqref{eq:wallcross}, the closed invariants, and the initial value $\ogw_{1,2}=2$.
\end{prop}
\begin{proof}
The closed case is given in~\cite[Claim 5.22]{KontsevichManin} and~\cite[Corollary 1]{RuanTian0}.
The open case is~\cite[Theorem 10]{ST3}.
\end{proof}

\section{The quadric -- preliminaries}\label{sec:quadtop}
\subsection{The singular cohomology}\label{ssec:dRquad}
\subsubsection{Notation}\label{sssec:notation}
For the rest of this paper, we compute the relative quantum cohomology and open Gromov-Witten invariants of $(X,L)$ where $X\subset \P^{n+1}$, $n>1$, is the quadric hypersurface and $L$ is as follows. Keeping in mind that all quadric hypersurfaces are isomorphic over $\C,$ we can take $X$ to be given by $\sum_{j=0}^{n}z_j^2 - z_{n+1}^2=0.$ We take $L$ to be the fixed locus of complex conjugation $\vec{z}\mapsto \vec{\bar{z}}.$ In particular, $L \simeq S^{n}.$ We consider the symplectic form $\omega$ given by the restriction of the Fubini-Study form $\oFS$ from $\C P^{n+1}.$
Let  $r:X\hookrightarrow \P^{n+1}$ be the inclusion and let
\[
r^*:H^*(\P^{n+1};\Z)\to H^*(X;\Z)
\]
be the induced restriction. Let $\check h \in H^*(\P^n;\Z)$ denote the Poincar\'e dual to a hyperplane $\hp \subset \C P^n.$ In other words, $\check h = [\oFS].$ Let
\[
h = r^* \check h.
\]
Let $p \in H^n(X;\Z)$ denote the Poincar\'e dual of $L.$

Let $m =\lceil n/2 \rceil$ and let $m^* = \lfloor n/2 \rfloor,$ so $m + \md = n.$ Let $\ml \subset X$ be a linear subspace of the maximal possible dimension. For concreteness, choose homogeneous coordinates $w_0, \ldots, w_{n+1},$ on $\C P^{n+1}$ such that $X$ is given by the equation
\begin{gather}
w_0w_1 + \cdots + w_{n-1}w_{n} + w_{n+1}^2 = 0, \qquad \text{when } n \text{ is odd, and} \notag %\label{eq:wodd}
\\
w_0w_1 + \cdots + w_{n}w_{n+1} = 0, \qquad \text{when } n \text{ is even.} \label{eq:weven}
\end{gather}
Take
\begin{equation*}%\label{eq:mlodd}
\ml   =  \{w_{2j} = 0 \,|\, j = 0,\ldots,m\} \subset \C P^{n+1}. \label{eq:mleven}
\end{equation*}
Let $l \in H^{2m}(X;\Z)$ be the Poincar\'e dual of~$\ml.$ When $n$ is odd, we will also use the linear subspace $\mla\subset X$ of maximal possible dimension given by
\begin{equation}\label{eq:mla}
\mla = \{w_{2j} = 0 \,|\, j = 0,\ldots,m - 1\} \cap \{w_{n+1} = 0\} \subset \C P^{n+1},
\end{equation}
and the Poincar\'e dual cohomology class $\la \in H^{2m}(X;\Z).$

\subsubsection{Cohomology ring}
It is well known that the singular cohomology of $X$ is described as follows.
\begin{prop}\label{prop:coh}
Equip the quadric hypersurface $X$ with its complex orientation. The ring structure of $H^*(X;\Z)$ is determined as follows.
For $n$ odd,
\[
H^*(X;\Z) \simeq \Z[h,l]/(h^m - 2l, l^2).
\]
For $n$ and $m$ even,
\[
H^*(X;\Z) \simeq \Z[h,l]/(h^{m+1} - 2hl, l^2 - h^ml).
\]
For $n$ even and $m$ odd,
\[
H^*(X;\Z) \simeq \Z[h,l]/(h^{m+1} - 2hl, l^2).
\]
For all $n,$ the Poincar\'e pairing of $H^*(X;\Z)$ is determined by $\int_X h^{\md} l = 1.$
\end{prop}

\begin{prop}\label{prop:quadmod}
As a group, the cohomology of $X$ is given by
\[
H^*(X;\Z)=\oplus_{j=0}^{2n}H^j(X;\Z),\qquad
H^j(X;\Z) \simeq\begin{cases}
0, & j\ne n, \text{ $j$ is odd},\\
\Z, & j\ne n, \text{ $j$ is even},\\
0, & j=n, \text{ $n$ is odd},\\
\Z\oplus \Z, & j=n, \text{ $n$ is even}.
\end{cases}
\]
It is generated by $\{h^j\}_{j=0}^{\md} \cup \{h^j l\}_{j = 0}^{\md}.$
\end{prop}

\begin{lm}\label{lm:ip}
We have
\begin{enumerate}
\item\label{it:bez}
$\int_X h^n = 2.$
\item\label{it:lin}
$\int_X h^{\md} l = 1.$
\item \label{it:del}
When $n$ or $m$ is odd,
$l^2 = 0.$ When $n$ and $m$ are even, $\int_X l^2 = 1.$
\end{enumerate}
\end{lm}
\begin{proof}
\begin{enumerate}
\item
By Bezout's theorem $\int_X h^n = \hp^n \cdot X = 2.$
\item
Since a codimension $\md$ linear subspace of $\ml$ is a single point, we have
\[
\int_X h^{\md}l =  \hp^{\md} \cdot \ml = 1.
\]
\item
First, consider the case that $n$ is odd, so $m > n/2.$ Consequently,
\[
|l^2| = 4m > 2n = \dim_\R X
\]
and $l^2 = 0.$

Now consider the case $n$ is even, so $m = n/2.$ Choose homogeneous coordinates $w_0,\ldots,w_{n+1},$ on $\C P^{n+1}$ such that $X$ is given by equation~\eqref{eq:weven} and let $\ml$ be as in~\eqref{eq:mleven}. Consider the maximal linear subspace $\ml' \subset X$ given by
\begin{gather*}
\ml' = \{w_{2j+1}= 0 \,|\, j = 0,\ldots, m\} \subset \C P^{n+1}, \qquad \text{when } m \text{ is odd, and} \\
\ml' = \{w_{n} = 0,\, w_{2j + 1} = 0 \,|\, j = 0,\ldots,m-1\} \subset \C P^{n+1}, \qquad \text{when } m \text{ is  even.}
\end{gather*}
We claim that $\ml$ and $\ml'$ are homologous in $X.$ Indeed, we can define a family of maximal linear subspaces $\ml_t \subset X, \; t \in [0,1],$ by
\[
\qquad \ml_t = \{t w_{3} - (1-t) w_0 = 0 = (1-t) w_{2} + t w_1, \, w_{2j} = 0\, |\,j = 2,\ldots,m\}.
\]
Then $\ml_0 = \ml$ and
\[
\ml_1 = \{w_1 = w_3 = w_4 = \cdots = w_{n} = 0 \} \subset \C P^{n+1}.
\]
Proceeding similarly with $w_k,\ldots,w_{k+3}$ instead of $w_0,\ldots,w_3,$ for $k = 4,8,\ldots$ successively, we conclude that $\ml$ and $\ml'$ are homologous as desired. Observe that $\ml \cap \ml' = \emptyset$ when $m$ is odd, but $\ml$ intersects $\ml'$ transversally at a single point when $m$ is even. Thus,
\[
\ml \cdot \ml = \ml \cdot \ml' =
\begin{cases}
0, & m \text{ odd,} \\
1, & m \text{ even.}
\end{cases}
\]
The desired formula for $l^2$ follows.
\end{enumerate}
\end{proof}

\subsubsection{Involution and maximal linear subspaces}
Let $\phi : X \to X$ be the involution given by
\begin{equation}\label{eq:hinv}
\phi([z_0:\cdots : z_{n+1}]) = [z_0: \cdots : -z_{n+1}].
\end{equation}
Recall the definition of $\mla$ and $\la$ from Section~\ref{sssec:notation}.
\begin{lm}\label{lm:phil}
For $n$ even, we have $\phi^* l = \la.$
\end{lm}
\begin{proof}
In the coordinates $w_i$ of equation~\eqref{eq:weven}, we have
\[
\phi([w_0:\cdots : w_n : w_{n+1}]) = [w_0: \cdots : w_{n+1}:w_n].
\]
So, $\phi(\ml) = \mla.$
\end{proof}

\begin{lm}\label{lm:lla}
For $n$ even, we have $h^m = l + \la \in H^{2m}(X).$
\end{lm}
\begin{proof}
The class $h^m$ is Poincar\'e dual to the intersection
\[
\{w_0 = 0\} \cap \{w_2 = 0\} \cap \cdots \cap \{w_{n-2} = 0\} \cap X = \ml + \mla.
\]
\end{proof}

\subsubsection{The fundamental class of the Lagrangian}
\begin{lm}\label{lm:lhs}
Let $W$ be a symplectic manifold of dimension $4w,$ let $K \subset W$ be a Lagrangian homology sphere, and let $k \in H^{2w}(W)$ be the Poincar\'e dual of $K.$ Then $\int_W k^2 =(-1)^{w}2.$
\end{lm}
\begin{proof}
Fix an orientation on $K$ and thus on $TK \simeq T^*K.$ Equip $TW$ with the canonical symplectic orientation. Let $N_K$ denote the normal bundle of $K \subset W.$ The short exact sequence
\[
0 \to TK \to TW \to N_K \to 0
\]
induces an orientation on $N_K.$ Explicitly, choose an oriented basis $u_1,\ldots,u_{2w} \in T_xK$
and choose $v_1,\ldots,v_{2w} \in T_xW$ such that $u_1,\ldots,u_{2w},v_1,\ldots,v_{2w}$ form a symplectic basis. That is,
\[
\omega_p(u_i,u_j) = 0 = \omega(v_i,v_j), \qquad \omega(u_i,v_j) = \delta_{ij}.
\]
Then $u_1 \wedge v_1 \wedge \ldots \wedge u_{2w} \wedge v_{2w}$ is the symplectic orientation of $T_xW.$ Let $\bar v_1,\ldots,\bar v_{2w},$ be the images of $v_1,\ldots,v_{2w},$ in $(N_K)_x.$ Then $(-1)^w \bar v_1 \wedge \ldots \wedge \bar v_{2w}$ is the induced orientation of $(N_K)_x.$ Indeed,
\[
u_1 \wedge v_1 \wedge \ldots \wedge u_{2w} \wedge v_{2w} = (-1)^{\frac{2w(2w-1)}{2}}u_1 \wedge \ldots \wedge u_{2w} \wedge v_1 \wedge \ldots \wedge v_{2w},
\]
and $(-1)^{\frac{2w(2w-1)}{2}}=(-1)^{w}.$

The symplectic form determines a bundle isomorphism $N_K \to T^*K.$ Explicitly, with $u^*_1,\ldots,u^*_{2w},$ the dual basis of $u_1,\ldots,u_{2w},$ the isomorphism is given by $\bar v_i \mapsto u_i^*.$ In particular, the isomorphism has sign $(-1)^w.$ Therefore,
\begin{multline*}
\int_W k^2 = [K]\cdot[K] = \int_K e(N_K)
= \\
=(-1)^{w}\int_K e(T^*K) = (-1)^{w}\chi(K) = (-1)^{w}\chi(S^{2w}) =(-1)^{w}2.
\end{multline*}
\end{proof}

Let $p \in H^n(X;\Z)$ denote the Poincar\'e dual of $L.$
\begin{lm}\label{lm:fcL}
When $n$ is odd, $p = 0.$ When $n$ is even, $p = \pm(h^m - 2l) = \mp (l - \la)$ with the sign depending on the orientation of $L.$
\end{lm}
\begin{proof}
When $n$ is odd, Proposition~\ref{prop:quadmod} gives $H^n(X) = 0,$ so $p = 0.$ When $n$ is even, since $L$ is Lagrangian, we know that
\begin{equation}\label{eq:yhm}
\int_X p \cup h^m = \int_X p \cup [\omega]^m = \int_L \omega^m = 0.
\end{equation}
By Lemma~\ref{lm:lhs}, we have
\begin{equation}\label{eq:y2}
\int_X p^2 = (-1)^{n/2} 2.
\end{equation}

By Proposition~\ref{prop:quadmod}, $\rk H^n(X;\Z) = 2.$ So, equations~\eqref{eq:yhm} and~\eqref{eq:y2} suffice to determine $p$ up to sign. Using Lemma~\ref{lm:ip}, one verifies these equations for $p = \pm(h^m - 2l).$ The equality $\pm(h^m - 2l) = \mp (l - \la)$ follows from Lemma~\ref{lm:lla}.
\end{proof}

\begin{prop}\label{pr:cohdrlag}
The de Rham cohomology of $X$ is given by
\[
H^*(X)\simeq \altfrac{\R[h]}{(h^{n+1})}
\]
if $n$ is odd
and
\[
H^*(X)\simeq \altfrac{\R[h,p]}{(h^{n+1},\quad hp,\quad p^2-(-1)^{n/2}h^n)}
\]
if $n$ is even. The Poincar\'e pairing is determined by $\int_X h^n = 2.$
\end{prop}

\begin{proof}
This follows from Proposition~\ref{prop:coh} and Lemma~\ref{lm:fcL}.
\end{proof}

\subsubsection{Relative cohomology ring}

Let $\At^*(L)$ be the reduced de Rham complex on $L$ and let $\Ht^*(L)$ the associated reduced cohomology; specifically, for $j>0$ we take $\At^j(L)=A^j(L)$, and for $j=0$ we take $\At^j(L)=A^0(L)/\R$, where $\R$ is thought of as the subspace of constant zero-forms. Define the associated relative complex $\At^*(X,L)$ by the short exact sequence
\[
0\rarr \At^*(X,L)\lrarr A^*(X)\lrarr \At^*(L) \rarr 0.
\]
Its cohomology $\Ht^*(X,L)$ then fits in the long exact sequence
\begin{equation}\label{eq:abc}
\xymatrix{
\widetilde H^*(X,L) \ar[rr]^a && H^*(X) \ar[dl]^{b} \\
 & \widetilde H^*(L) \ar[ul]_{[1]}^c.
}
\end{equation}
Let $\tilde u \in \widetilde H^{n+1}(X,L)$ be the image under $c$ of the class with integral $1$ in $\widetilde H^n(L).$ Let $\tilde h = [r^*\omega_{FS}] \in \widetilde H^2(X,L).$
\begin{prop}\label{pr:Ht}
As $\R$ algebras, when $n$ is odd
\[
\widetilde H^*(X,L) \simeq \R[\tilde h,\tilde u]/(\tilde h^{n+1},\tilde u^2,\tilde u\tilde h).
\]
When $n$ is even,
\[
\widetilde H^*(X,L) \simeq \R[\tilde h]/(\tilde h^{n+1}).
\]
\end{prop}
\begin{proof}
By definition, we have $a(\tilde h) = h.$ Since $L \simeq S^n,$ it follows that $\widetilde H^*(L)$ is one dimensional and concentrated in degree $n.$
For the case $n$ odd, observe that $b = 0$ for degree reasons, so $c$ is injective and $a$ is an isomorphism except in degree $n+1.$ By exactness, $\tilde{u}\in \Im(c)$ implies $a(\tilde{u})=0$. Since $a(\tilde u\tilde h) = a(\tilde u)a(\tilde h) = 0,$ we deduce that $\tilde u\tilde h = 0.$ Also, $\tilde u^2 = 0$ for degree reasons. The claim follows from Proposition~\ref{pr:cohdrlag}.

For the case $n$ even, the integral of $b(p) \in \widetilde H^*(L)$ is $(-1)^{n/2}2$ by Lemma~\ref{lm:lhs}. In particular, $b(p) \neq 0,$ so $b$ is surjective and $c = 0.$ So, the long exact sequence~\eqref{eq:abc} implies that $\widetilde H^*(X,L)$ is the kernel of~$b,$ which is the subalgebra of $\widetilde H^*(X)$ generated by $\tilde h.$ So, the claim follows by Proposition~\ref{pr:cohdrlag}.
\end{proof}

Recalling the long exact sequence~\eqref{eq:yrho}, let $\hat u = y(1) \in \Hh^{n+1}(X).$ Let $\hat h = [r^*\omega_{FS}] \in \Hh^2(X,L).$
\begin{prop}\label{pr:Hh}
As $\R$-vector spaces, when $n$ is odd,
\[
\Hh^*(X,L) \simeq \R[\hat h,\hat u]/(\hat h^{n+1},\hat u^2, \hat u \hat h).
\]
When $n$ is even,
\[
\Hh^*(X,L) \simeq \altfrac{\R[\hat h]}{(\hat h^{n+1})}.
\]
\end{prop}
\begin{proof}
Since $L = S^n,$ Lemma~5.14 of~\cite{ST2} asserts that the inclusion of underlying complexes induces an isomorphism $\widetilde H^*(X,L) \overset{\sim}{\to} \Hh^*(X,L).$ So, the claim follows from Proposition~\ref{pr:Ht}.
\end{proof}

\begin{rem}
Proposition~\ref{pr:Hh} gives a natural splitting of the exact sequence~\eqref{eq:splitL0}. Indeed, a section of $\rho$ is given by mapping $h^j \in H^*(X)$ to $\hat h^j \in \Hh^*(X)$ for $j = 0,\ldots,n.$
\end{rem}

\subsection{Additional topological properties}\label{ssec:top}

\begin{lm}\label{lm:varpi}
We have $H_2(X,L;\Z)\simeq \Z$. Moreover, the map $\varpi:H_2(X;\Z)\to H_2(X,L;\Z)$ from the long exact sequence of the pair $(X,L)$ is an isomorphism if $n>2$ and a surjection if $n=2$.
\end{lm}

Identify $H_2(X,L;\Z)$ with $\Z$ so that the non-negative integers correspond to classes $\beta \in H_2(X,L;\Z)$ with $\omega(\beta) \geq 0$. When $n>2$, we also identify the non-negative integers with classes $\hat\beta\in H_2(X;\Z)$ with $\omega(\beta)>0$. Then $\varpi=\Id_\Z$. When $n=2$, Lemma~\ref{lm:fcL} implies that the inclusion $H_2(L;\Z)\hookrightarrow H_2(X;\Z)$ is given by $a\mapsto \pm(a,-2a)$. Exactness then gives $\varpi(a,b)=2a+b$.

\begin{rem}
Using the identifications above, we can write explicitly the integrals of divisors, analogously to equation~\eqref{eq:divpn}. Namely,
by Proposition~\ref{pr:Hh}, for any $\beta\in H_2(X,L;\Z)$ we have
\[
\int_\beta \hat{h}=\beta.
\]
By Proposition~\ref{prop:quadmod}, when $n>2$ we have
\[
\int_{\hat\beta} h =\hat\beta,\qquad \forall\hat\beta\in H_2(X;\Z),
\]
and when $n=2$ we have
\[
\int_{(\beta_1,\beta_2)}h=\beta_1,\quad \int_{(\beta_1,\beta_2)}l=\beta_2, \qquad \forall (\beta_1,\beta_2)\in H_2(X;\Z).
\]
As in Section~\ref{sec:cpn}, these values will be used when applying the divisor axioms~\eqref{ax_divisor} and~\eqref{ax_c_divisor}.
\end{rem}

\begin{lm}\label{lm:char}
The relevant characteristic classes of $(X,L)$ are as follows:
\begin{enumerate}
\item\label{it:c1}
For all $\hat\beta\in H_2(X;\Z)$ we have
$
c_1(\hat\beta)=n\hat\beta.
$
\item\label{it:mu}
For all $\beta\in H_2(X,L;\Z)$, we have
$
\mu(\beta)=2n\beta.
$
\item\label{it:ws}
$w_\s=0$.
\end{enumerate}
\end{lm}

\subsection{Basic values for GW and OGW invariants}\label{ssec:quadbase}
Here we compute the basic values required to deduce the quantum ring structure on the absolute and relative cohomologies.
In this section we assume $n\ge 3$. For the case $n=2$ see Section~\ref{sec:quad2}. In particular, we identify $H_2(X;\Z)$ and $H_2(X,L;\Z)$ with $\Z$, so that positive-energy classes correspond to $\Z_{>0}$.

For $\beta=0$, the invariant $\GW_\beta(h^i,h^j,h^m)$ does not vanish if and only if $i+j+m=n$, in which case by Lemma~\ref{lm:ip}\ref{it:bez}
we have $\GW_0(h^i,h^j,h^m)=2$.

\begin{lm}\label{lm:gwnn-1}
$\GW_1(h^{m^*}l,h^{m^*-1}l)=1$.
\end{lm}

\begin{proof}
By Lemma~\ref{lm:ip}\ref{it:lin}, the question translates to the number of lines in $X$ though one point and one line.
Any point $x\in X$ and line $\ell\subset X$ determine a plane in $\P^{n+1}$. This plane intersects the hypersurface $X$ in a degree $2$ curve $C$. On the other hand, $C$ contains the line $\ell$. Therefore, $C$ is the union of two lines, $\ell\cup \ell'$. In other words, there is a unique line $\ell'$ in $X$ passing through $x$ that intersects $\ell$. Thus, $\GW_1(h^{m^*}l, h^{m^*-1}l)=1$.
\end{proof}

\begin{cor}\label{cor:val}
$\GW_1(h^n, h^{n-1})=2^2\cdot \GW_1(h^{m^*l},h^{m^*-1}l)=4. $
\end{cor}

\begin{lm}\label{lm:gwd1}
For $i,j,k,$ such that $i+j+k=2n$, $i\le j\le k$, we have
\[
\GW_1(h^i,h^j,h^k)=\begin{cases}
0, & (i,j,k)=(0,n,n),\\
4, & i>0, k=n,\\
8, & i>0, k<n.
\end{cases}
\]
\end{lm}
Note that the condition $i+j+k=2n$ is required by the degree axiom~\eqref{ax_c_deg}.

\begin{proof}

For $i=0$,  we apply the unit axiom~\eqref{ax_c_unit} to get the first case:
\begin{equation*}%\label{eq:gw10}
\GW_1(h^0,h^n,h^n)=0.
\end{equation*}
For $i=1$, we apply the divisor axiom~\eqref{ax_c_divisor} and Corollary~\ref{cor:val} to get
\begin{equation*}%\label{eq:gw14}
\GW_1(h,h^{n-1},h^n)=\GW_1(h^{n-1},h^n)=4.
\end{equation*}

Now assume $i>1$.
Take the coefficient of $T^1$ in the equation~\eqref{eq:wdvv} with $u=h^{i-1}, v=h, w=h^j,$ and $y=h^k,$ to get
\[
\GW_1(h^{i},h^{j},h^{k})+\GW_1(h^{i-1},h^{j+k}) =
\GW_1(h^{i-1}, h^{j+1}, h^{k})+ \GW_1(h^{j},h^{i+k-1}).
\]

The summand $\GW_1(h^{i-1},h^{j+k})$ contributes nothing:
By the degree axiom~\eqref{ax_c_deg}, it can only be nonzero if $j+k=n$ and $i-1=n-1$, but in this case $i=n>j,$ contradicting our assumption.

The summand $\GW_1(h^{j},h^{i+k-1})$ satisfies
\[
\GW_1(h^{j},h^{i+k-1}) = \begin{cases}
4, & (i,j,k)= (2, n-1, n-1),\\
0, & otherwise,
\end{cases}
\]
for the following reason. By the degree axiom~\eqref{ax_c_deg}, $\GW_1(h^{j},h^{i+k-1})$ can only be nonzero if $\{j,i+k-1\}=\{n-1,n\}$. Since we assumed $i>1$, we have $i+k-1>k\ge j$, so $(j,i+k-1)=(n-1,n)$. In this case $k=n+1-i$. Since $i>1$ we get $k<n$, and since $k\ge j=n-1$, we get $k=n-1$ and $i=2$.

Thus,
if $(i,j,k)=(2,n-1,n-1)$, we have
\[
\GW_1(h^i,h^j,h^k)=
\GW_1(h^{2},h^{n-1},h^{n-1}) =
\GW_1(h, h^{n}, h^{n-1})+ \GW_1(h^{n-1},h^{n}) =8,
\]
and if $(i,j,k)\ne(2,n-1,n-1)$ with $i>1$, we have
\begin{equation}\label{eq:indc}
\GW_1(h^{i},h^{j},h^{k}) =
\GW_1(h^{i-1}, h^{j+1}, h^{k}).
\end{equation}
Thus, if $k=n$ then $i+j=n$, and inductively,
\[
\GW_1(h^i,h^{n-i},h^n)
= \GW_1 (h^{i-(i-1)},h^{(n-i)+(i-1)}, h^n)
= \GW_1 (h, h^{n-1}, h^n)=4.
\]
If $k<n$, we would again like to apply equation~\eqref{eq:indc} repeatedly. In order to do so,
consider the transformation
\[
\psi:\Z^{\oplus 3} \lrarr \Z^{\oplus 3}
\]
defined by
\[
\psi(a,b,c):=\begin{cases}
(a-1,b+1,c), & b+1\le c,\\
(a-1,c,b+1), & b+1>c.
\end{cases}
\]
Note that $\psi$ preserves the sum of the components, and takes non-decreasing triples to non-decreasing ones.
Since $k<n$, there exists $q\in\Z_{\ge 0}$ such that
$\psi^q(i,j,k)$ has the form $\psi^q(i,j,k)=(\nu,\mu,n-1)$.
Then $\mu=n+1-\nu$, $\mu\le n-1$, and consequently $\nu\ge 2$. Taking $q':=\nu-2$, we get
\[
\psi^{q+q'}(i,j,m) = \psi^{q'}(\nu,\mu,n-1)
=(\nu-(\nu-2),(n+1-\nu)+(\nu-2),n-1)
=(2,n-1,n-1).
\]
Thus,
\[
\GW_1(h^i,h^{j},h^k) =
\GW_1(h^2,h^{n-1},h^{n-1}) = 8.
\]
\end{proof}

\begin{lm}\label{lm:11n}
$\ogw_{1,1}(h^n)=2$.
\end{lm}

\begin{proof}
Apply $\d_s$ to equation~\eqref{eq:cor2} with $v=h$ and $w=h^{n-1}$, set $t_j=s=0$ and take the coefficient of $T^1$. This gives
\[
\frac{1}{2}\GW_0(h,h^{n-1},1)\ogw_{1,1}(h^n) + \frac{1}{2}\GW_1(h,h^{n-1},h^n)\ogw_{0,1}(1) = 0.
\]
Apply the zero axiom~\eqref{ax_c_zero} and Lemma~\ref{lm:ip}\ref{it:bez} to the first summand, and apply the divisor~\eqref{ax_c_divisor} and zero~\eqref{ax_zero} axioms and Corollary~\ref{cor:val} to the second summand, to get
\[
\ogw_{1,1}(h^n) - 2 = 0.
\]

\end{proof}

\section{The quadric -- odd dimension}

\subsection{Absolute quantum ring}\label{ssec:abquadodd}
In this section, we assume $n>1$ is odd. We work with the basis
\[
\dl_j:= h^j,\qquad j=0,\ldots,n,
\]
of $H^*(X;\R)$. By Lemma~\ref{lm:ip}\ref{it:bez}, the Poincar\'e pairing on $H^*(X;\R)$ is given by
\[
g_{ij}= \int_X \dl_i\wedge\dl_j
=\delta_{i,n-j}\cdot 2,
\]
where $\delta_{k,l}$ is the Kronecker delta, and the inverse matrix $(g^{ij})_{i,j}:={(g_{ij})_{i,j}}^{-1}$ is given by
$g^{ij}=\delta_{i,n-j}\cdot \frac{1}{2}$.
Thus, the quantum product is given by
\[
\dl_i*\dl_j =  \sum_{\substack{0\le m\le n\\ d\ge 0}} \frac{1}{2}T^{d}\GW_d(\dl_i,\dl_j,\dl_m)\dl_{n-m},
\qquad i,j\in\{0,\ldots,n\}.
\]
\begin{prop}\label{prop:qhquadod}
The absolute quantum cohomology is given by
\[
QH^*(X)\simeq \altfrac{\R[[q]][x]}{(x^{n+1}-4qx)}.
\]
\end{prop}
Here, as in the case of $\P^n$, the variable $q$ corresponds to $T^{\varpi(\hat\beta_1)}$ where $\hat\beta_1$ is the class of a line generating $H_2(X;\Z)$, and $x$ corresponds to a divisor generating $H^2(X)$. Note that for the quadric, $T^{\varpi(\hat\beta_1)}=T^1=T$.
\begin{proof}
If $i+j<n$, use the degree axiom~\eqref{ax_c_deg} and the zero axiom~\eqref{ax_c_zero} to get
\begin{align*}
\dl_i*\dl_j=&\sum_{d,m}\frac{1}{2}T^d\GW_d(\dl_i,\dl_j,\dl_m)\dl_{n-m}\\
=& \sum_m\frac{1}{2}\cdot \GW_0(\dl_i,\dl_j,\dl_m)\dl_{n-m}\\
=& \frac{1}{2}\cdot \big(\int_{X}\dl_i\wedge\dl_j\wedge\dl_{n-i-j}\big)\cdot \dl_{i+j}\\
=&\dl_{i+j}.
\end{align*}
By the associativity of $*$, for $k<n$ we therefore have
\[
\dl_1^{*k}=\dl_k.
\]
Next, we compute $\dl_1^{* n}=\dl_1^{* (n-1)}*\dl_1$, as follows. The degree axiom~\eqref{ax_c_deg} implies that $\GW_d(\dl_{n-1},\dl_1,\dl_m)$ can only have nonzero values for $d=0,1$. By the zero axiom~\eqref{ax_c_zero}, the case $d=0$ contributes $\GW_0(\dl_{n-1},\dl_1,\dl_0)=2$, as before. For $d=1$, again by the degree axiom~\eqref{ax_c_deg}, the only nonzero case is when $m=n$. By the divisor axiom~\eqref{ax_c_divisor}, the value is
\[
\GW_1(\dl_{n-1},\dl_1,\dl_n)
=\int_{\beta_1}\dl_1\cdot \GW_1(\dl_{n-1},\dl_n)
=\GW_1(\dl_{n-1},\dl_n)=4.
\]
So,
\[
\dl_1^{* n}
=\frac{1}{2}\GW_0(\dl_{n-1},\dl_1,\dl_0)\dl_n
+ \frac{1}{2}T\GW_1(\dl_{n-1},\dl_1,\dl_n)\dl_0
=\dl_n+2T.
\]
Thus, we see that $\{\dl_1^{* j}\}_{j=0}^n$ generate the $Q_X$-module $QH^*(X)$. To complete the computation of the ring structure, we need to compute powers of $\dl_1$ higher than $n$. By the degree axiom~\eqref{ax_c_deg} and the zero axiom~\eqref{ax_c_zero}, we have
\[
\dl_n*\dl_1
=\frac{1}{2}T\GW_1(\dl_n,\dl_1,\dl_{n-1})\dl_{1}
=\frac{1}{2}T\cdot 4\dl_{1} =2T\dl_1.
\]
Therefore,
\[
\dl_1^{*(n+1)}
= \dl_1^{* n}*\dl_1 = (\dl_n+2T)*\dl_1
= 2T\dl_1+2T\dl_1=4T\dl_1.
\]
Taking the variables $q,x,$ to correspond to $T,\dl_1,$ respectively, we get
the desired result.
\end{proof}

\subsection{Relative quantum ring}\label{ssec:relquadodd}

In this section, we work with the following basis.
Let $\g_1=[r^*\oFS]\in \Hh^*(X,L;\R)$ and set $\g_j=\g_1^{\wedge j}$ for $j=0,\ldots,n$. Let $\g_\diamond=y(1)$ where $y$ is the map from~\eqref{eq:yrho}.
Then $\{\g_j\}_{j=0}^n\cup\{\g_\diamond\}$ is a basis of $\Hh^*(X,L;\R)$, and for $\dl_j$ of Section~\ref{ssec:abquadodd} we have
\[
\rho(\g_j)=\dl_j,
\qquad j=0,\ldots,n.
\]
For convenience, we also define $\dl_\diamond:=0\in H^*(X;\R)$.

For $i,j\in \{0,\ldots,n\}\cup\{\diamond\}$, we have
\[
\mem(\g_i,\g_j)
= \sum_{\substack{0\le m\le n\\ d\ge 0}} \frac{1}{2}T^{d}\GW_d(\dl_i,\dl_j,\dl_m)\g_{n-m}
+ \sum_{d\ge 0}T^d\ogwb_{d,0}(\g_i,\g_j)\g_\diamond.
\]

\begin{proof}[Proof of Theorem~\ref{thm:rqhqodd}]
By the degree axiom~\eqref{ax_deg} of $\ogw$, the only nonzero values possible for $\ogw_{d,0}(\g_i,\g_j)$ are when $d=0,1$. When $d=0$, by the zero axiom~\eqref{ax_zero} of $\ogw$ and the wall-crossing theorem~\eqref{eq:wallcross}, the only nontrivial value is
\[
\ogw_{0,0}(\g_\diamond,\g_0)=-\ogw_{0,1}(\g_0)=1.
\]
When $d=1$, by the degree axiom~\eqref{ax_deg} nonzero values are only possible when $|\g_i|+|\g_j|=3n+1$. If, without loss of generality, $j=\diamond$, we use the wall-crossing theorem~\eqref{eq:wallcross} and Lemma~\ref{lm:11n} to compute
\[
\ogw_{1,0}(\g_n,\g_\diamond)
=-\ogw_{1,1}(\g_n)=-2.
\]
If $i,j\in \{0,\ldots,n\}$, take the coefficient of $T$ in the first OWDVV equation with $u=t_i$, $w=t_1$, and $v=t_{j-1}$:
\[
\ogw_{1,0}(\g_i,\g_j) = \ogw_{1,0}(\g_{i+1},\g_{j-1}).
\]
In the case $i=n$ we interpret the right-hand side by setting $\g_{n+1}:=0$, and the identity still holds.
Note that $2\cdot\big((i+1)+(j-1)\big)=2\cdot(i+j)=3n+1$. By induction, we increase the degree of the first constraint $i$ (while decreasing the second) until it reaches $n$
 and conclude
$\ogw_{1,0}(\g_i,\g_j)=0$.

Since $\rho$ is an isomorphism when restricted to $W=Span(\g_0,\ldots,\g_n)$, we can summarize
\[
\mem(\g_i,\g_j) = \begin{cases}
\g_\diamond, & \{i,j\}=\{\diamond, 0\},\\
-2T\g_\diamond, & \{i,j\}=\{\diamond,n\},\\
(\rho|_W)^{-1}(\dl_i*\dl_j), & otherwise.
\end{cases}
\]
Taking $q,x,y,$ to be the variables corresponding to $T,\g_1,\g_\diamond,$ respectively, we get
the desired result.
\end{proof}

\subsection{Full computation of invariants}
In this section we prove Theorem \ref{thm:oddcomputability}.

In Appendix~\ref{sec:app}, we compute the following value.
\begin{lm}\label{lm:init}
The unique spin structure on $L$ is such that $\ogw_{1,3}=(-1)^{\frac{(n-2)(n-3)}{2}}2$.
\end{lm}

Let $n = 2r+1$. From Lemma \ref{lm:gwd1} we have another initial value:
\begin{lm}
	\label{lm:singlegwoddquadric}
	$\GW_1(\dl_{r+1}, \dl_{r+1}, \dl_{2r}) = 8$.
\end{lm}

We are now ready to prove computability lemmas.
\begin{lm}
	\label{lm:oddogwsingleinterior}
	All open Gromov-Witten invariants $\ogw_{\beta,k}(\gamma_j)$ with $j \neq 1, \diamond$ can be computed from closed Gromov-Witten invariants and open Gromov-Witten invariants of strictly lower degree.
\end{lm}
\begin{proof}
	 Consider the invariant $\ogw_{\beta,k}(\Gamma_j)$. If $j = 0$, the invariant vanishes unless $\beta = 0$, in which case the value is determined by the $\ogw$ axioms. Thus we can assume $j\geq 2$. Note that for degree reasons we must have $k \geq 1$. We can then consider the open WDVV equation \eqref{eq:cor2} with $v=1$, $w = j-1$ in degree $\beta$ after applying $\partial_s^{k-1}$:
	\begin{align}
		& \GW_0(\dl_1, \dl_{j - 1}, \dl_{n-j})\ogw_{\beta, k}(\Gamma_{j})g^{n-j,j} \notag\\
		\label{1 interior constraint term 1}
		& + \sum_{\substack{\beta_1 + \beta_2 = \beta \\ \beta_2 < \beta\\ i = 0, \dots, n}} \GW_{\beta_1}(\dl_1, \dl_{j - 1},\dl_i)\ogw_{\beta_2, k}(\Gamma_{n-i})g^{n-i,i}  \\
		\label{1 interior constraint term 2}
		& - \sum_{\substack{\beta_1 + \beta_2 = \beta\\ \beta_1, \beta_2 < \beta \\ i = 1, \dots, k-1}} {k-1 \choose i} \ogw_{\beta_1, 1+k-i}\ogw_{\beta_2,i}(\Gamma_1, \Gamma_{j - 1})  \\
		\label{1 interior constraint term 3}
		= & - \sum_{\substack{\beta_1 + \beta_2 = \beta\\  \beta_1, \beta_2 < \beta \\ i = 0, \dots, k-1}} {k-1 \choose i} \ogw_{\beta_1, 1+i}(\Gamma_1)\ogw_{\beta_2,k-i}(\Gamma_{j - 1})
	\end{align}
	Here in line \eqref{1 interior constraint term 1}, we have $\beta_2 < \beta$, as the only (possibly) non-zero term with $\beta_2 = \beta$ involves the term we are trying to compute. In line \eqref{1 interior constraint term 2} both $\beta_1, \beta_2 < \beta$, as there are no non-trivial $\ogw_{0,k}$ with either 0 or 2 interior constraints, which follows from the choice of $P_{\mathbb{R}}$ and the degree-zero property. Finally in line \eqref{1 interior constraint term 3}, we have  $\beta_1, \beta_2 < \beta$ as the the invariant with $\beta_2 = 0$ that could appear $\ogw_{0,k-i}(\Gamma_{j-1})$ vanishes by the energy-zero property, as $j > 1$.
	As $GW_0(\dl_1, \dl_{j - 1}, \dl_{n-j}) = 2$, we have expressed $\ogw_{\beta, k}(\Gamma_j)$ in terms of closed Gromov-Witten invariants and open Gromov-Witten invariants of strictly lower degree.
\end{proof}
\begin{lm}
	\label{lm:oddogwmultipleinterior}
	An open Gromov-Witten invariant $\ogw_{\beta,k}(\Gamma_{j_1}, \dots , \Gamma_{j_l})$ with $l\geq 2$ can be computed from closed Gromov-Witten invariants, open Gromov-Witten invariants of smaller degree, and invariants of the form $\ogw_{\beta,k'}(\Gamma_{j_1 - 1},\Gamma_J)$ for $k' \leq k$ and some index set $J$.
\end{lm}
\begin{proof}
	By reordering the inputs, we may assume $j_1 \leq j_2 \leq \dots \leq j_l$. Consider the open WDVV equation \eqref{eq:cor1} with $u = j_2$, $v = j_1 - 1$, $w = 1$, and take the coefficient of $T^\beta t_{j_3} \dots t_{j_l}$ after applying $\partial_s^k$. Set $I = (j_3, \dots, j_l)$. This yields:
	\begin{align}
		\label{many interior constraints term 0}
		& \GW_0(\Gamma_1, \Gamma_{j_1 - 1}, \Gamma_{n-j_1})\ogw_{\beta, k}(\Gamma_{j_1}, \dots, \Gamma_{j_l})g^{n-j_1,j_1} \\
		\label{many interior constraints term 1}
		& + \sum_{\substack{\beta_1 + \beta_2 = \beta \\ \beta_2 < \beta \\ I_1 \sqcup I_2 = I \\ i = 0, \dots, n}} \GW_{\beta_1}(\Gamma_1, \Gamma_{j_1 - 1}, \Gamma_{I_1}, \Gamma_i)\ogw_{\beta_2, k}(\Gamma_{n-i}, \Gamma_{j_2}, \Gamma_{I_2})g^{n-i,i} \\
		\label{many interior constraints term 2}
		& - \sum_{\substack{\beta_1 + \beta_2 = \beta\\  \beta_1, \beta_2 < \beta \\ I_1 \sqcup I_2 = I \\ i = 0, \dots, k}} {k \choose i} \ogw_{\beta_1, 1+i}(\Gamma_{j_2}, \Gamma_{I_1})\ogw_{\beta_2,k-i}(\Gamma_{j_1 - 1}, \Gamma_1, \Gamma_{I_2}) \\
		\label{many interior constraints term 5}
		& = \GW_0(\Gamma_{j_2}, \Gamma_{1}, \Gamma_{n-j_2 - 1})\ogw_{\beta, k}(\Gamma_{j_2 + 1}, \Gamma_{j_1 - 1}, \Gamma_I)g^{n-j_2-1,j_2+1}\\
		\label{many interior constraints term 3}
		& + \sum_{\substack{\beta_1 + \beta_2 = \beta \\ \beta_2 < \beta \\ I_1 \sqcup I_2 = I \\ i = 0, \dots, n}} \GW_{\beta_1}(\Gamma_{j_2}, \Gamma_{1}, \Gamma_{I_1}, \Gamma_i)\ogw_{\beta_2, k}(\Gamma_{n-i}, \Gamma_{j_1 - 1}, \Gamma_{I_2})g^{n-i,i}\\
		\label{many interior constraints term 4}
		&- \sum_{\substack{\beta_1 + \beta_2 = \beta\\  \beta_1, \beta_2 < \beta \\ I_1 \sqcup I_2 = I \\ i = 0 , \dots, k}} {k \choose i} \ogw_{\beta_1, 1+i}(\Gamma_{j_1 - 1} \Gamma_{I_1})\ogw_{\beta_2,k-i}(\Gamma_{j_2}, \Gamma_1, \Gamma_{I_2})
	\end{align}
	Here in line \eqref{many interior constraints term 1}, we have $\beta_2 < \beta$, as the only (possibly) non-zero term with $\beta_2 = \beta$ involves the term we are trying to compute, line \eqref{many interior constraints term 0}. For lines \eqref{many interior constraints term 2} and \eqref{many interior constraints term 4}, we need $\beta_1, \beta_2 > 0$, as the $\ogw$ have at least one interior constraint on $\Gamma_{j_1 - 1}$ or $\Gamma_{j_2}$, and by assumption $j_1, j_2 > 1$, so that by the degree zero axiom, there are no contributions from these terms in degree $0$. This implies that $\beta_1,\beta_2 < \beta$ for lines \eqref{many interior constraints term 2} and \eqref{many interior constraints term 4}. Line \eqref{many interior constraints term 3} does not contain terms with $\beta_2 = \beta$, as these appear in the line we have singled out \eqref{many interior constraints term 5}
	We have thus reduced the invariant $\ogw_{\beta,k}(\Gamma_{j_1}, \dots , \Gamma_{j_l})$ to terms of lower degree, or of the same degree, but with the degree of one of the interior constraints decreased by 2.
\end{proof}
Using the above two lemmas we can compute:
\begin{lm}
	\label{lm:singleogwoddquadric}
	The closed Gromov-Witten invariants, together with Lemmas \ref{lm:oddogwsingleinterior} and \ref{lm:oddogwmultipleinterior} allow us to compute: $\ogw_{1,1}(\Gamma_{r+1}, \Gamma_{r+1}) = 2$.
\end{lm}
\begin{lm}
	\label{lm:oddogwnointerior}
	All open Gromov-Witten invariants with no interior constraints, except for $\ogw_{1,3}$, can be computed from invariants of lower degree, and invariants of equal degree, but with at least one interior constraint of degree greater than 2.
\end{lm}
\begin{proof}
	Consider the open WDVV equation \eqref{eq:cor2}, with $v = w = r+1$, and apply $\partial_s^{k-1}$, taking the coefficient of $T^{\beta + 1}$. We find:
	\begin{align}
		\label{no interior constraints term 0}
		& \GW_1(\dl_{r+1}, \dl_{r+1}, \dl_{2r})\ogw_{\beta, k}(\Gamma_1)g^{2r,1} \\
		\label{no interior constraints term 1}
		& + \sum_{\substack{\beta_1 + \beta_2 = \beta+1 \\ \beta_2 < \beta\\ i = 0, \dots, n}} GW_{\beta_1}(\Gamma_{r+1}, \Gamma_{r+1},\Gamma_i)\ogw_{\beta_2, k}(\Gamma_{n-i})g^{i,n-i}\\
		\label{no interior constraints term 4}
		& - (k-1)\ogw_{\beta,k}\ogw_{1,1}(\Gamma_{r+1}, \Gamma_{r+1})\\
		\label{no interior constraints term 2}
		& - \sum_{\substack{\beta_1 + \beta_2 = \beta+1 \\\beta_1 < \beta\\\beta_2 \leq \beta \\ i = 1, \dots, k-1}} {k-1 \choose i} \ogw_{\beta_1, 1+k-i}\ogw_{\beta_2,i}(\Gamma_{r+1},\Gamma_{r+1})\\
		\label{no interior constraints term 3}
		= & - \sum_{\substack{\beta_1 + \beta_2 = \beta + 1\\ \beta_1, \beta_2 \leq \beta \\ i = 0, \dots, k-1}} {k-1 \choose i} \ogw_{\beta_1, 1+i}(\Gamma_{r+1})\ogw_{\beta_2,k-i}(\Gamma_{r+1})
	\end{align}
	For line \eqref{no interior constraints term 1}, note there are no contributions with $\beta_2 = \beta + 1$, as $GW_0(\Gamma_{r+1},\Gamma_{r+1},\Gamma_i) = 0$ for any $i$, by the degree property. There are also no contributions with $\beta_2 = \beta$, as the only $GW_1(\Gamma_{r+1},\Gamma_{r+1},\Gamma_i)$ that contributes is the one with $i = 2r$, which is the term in the line above it.
Secondly for line \eqref{no interior constraints term 2} there are no contributions with $\beta_1$ or $\beta_2 = \beta + 1$, as then the other term is $\ogw_{0, \star}$ with either 0 or 2 interior marked points, which thus vanishes by the degree zero property. Also, there only contribution with $\beta_1 = \beta$ has been singled out in the line above it.
Finally for \eqref{no interior constraints term 3} there are no contributions with $\beta_1$ or $\beta_2 = \beta + 1$, as then the other term is $\ogw_{0, \star}$ with 1 interior marked point on $\Gamma_{r+1}$. As $r+1 > 0$, this vanishes by the degree zero property. By the previous lemma, we can combine lines \eqref{no interior constraints term 0} and \eqref{no interior constraints term 4} to: $2(2\beta - k+1)\ogw_{\beta,k}$.

By the degree property of the $\ogw$, we have $2n\beta = k(n-1) + 3 - n$, and thus $n(2\beta - k + 1) = 3-k \neq 0$ unless $k = 3$ and $\beta = 1$. The invariant $\ogw_{1,3}$ is precisely the invariant we have computed already, thus the above equations determine $\ogw_{\beta, k}$ in terms of invariants of lower degree, or of equal degree, with at least 1 interior constraint.
\end{proof}

Finally, we combine the above information to prove Theorem \ref{thm:oddcomputability}.
\begin{proof}[Proof of theorem \ref{thm:oddcomputability}]
	Let $\ogw_{\beta, k}(\Gamma_I)$ be an invariant we want to compute. By the divisor axiom, we may assume $1 \notin I$. By Theorem \ref{thm:wallcross} we may assume $\diamond \notin I$. Furthermore, the number of boundary constraints $k$ is determined by $\beta, I$ and the degree axiom of $\ogw$. We now introduce a total ordering on the set $\{(\beta, I) \in \mathbb{N} \times \{0,2,3,\dots, n\}^{|I|}$ of parameters for open Gromov-Witten invariants.
	Order the degrees $\beta$ by the natural ordering on $\mathbb{N}$ (under the identification $H_2(X,L) \cong \mathbb{Z}$). Order the index sets $I$ as follows: the element $I = \emptyset$ is the largest element. Then, if $|I| \neq \emptyset$ order such $I$ by $|I|$ under the natural ordering on $\mathbb{N}_{>0}$. Assume that $I$ is ordered so that $i_1 \leq i_2 \leq \dots \leq i_{|I|}$. If two sets $I, I'$ have the same size, say $I < I'$ if $I$ appears first in the lexicographic ordering on $(\mathbb{N} \setminus \{1\})^{|I|}$. The total order on tuples $(\beta, I)$ is then given by the lexicographic order. This total order satisfies the well-ordering principle, so we can apply induction. The three lemmas \ref{lm:oddogwsingleinterior}, \ref{lm:oddogwmultipleinterior} and \ref{lm:oddogwnointerior} precisely show that the open WDVV equations allow us to compute the invariant $\ogw_{\beta, k}(\Gamma_I)$ in terms of invariants strictly lower in the ordering on $(\beta, I, k)$. Except for the case $\ogw_{1,3}$ which we calculated separately. The result thus follows by induction.
\end{proof}
We now consider how the invariant change if instead we run the above induction process with $OGW_{1,3} = -2$. We first have the following observation:
\begin{lm}
	\label{lm: even boundary constraints vanish}
	For $I \subset \{0,1,\dots, n\}^{l}$ and $k$ even, we have: $OGW_{\beta, k}(\Gamma_I) = 0$.
\end{lm}
\begin{proof}
	The result holds for $\beta = 0$. Then apply the above recursion process, noting that always at least one of the terms in the products has an even number of boundary constraints.
\end{proof}
We then show: \begin{prop}
	\label{prop: sign change odd quadrics}
	Suppose that instead we run the above induction process with the initial value $\ogwb_{1,3} = -2$, then for any $I \subset \{0,1,\dots, n\}^{l}$, the resulting invariant $\ogwb_{\beta,k}(\Gamma_I)$ changes by the sign $(-1)^{\frac{k-1}{2}}$ and the invariant $\ogw_{\beta,k}(\Gamma_\diamond^{\otimes r},\Gamma_I)$ changes by the sign $(-1)^{\frac{k+r-1}{2}}$.
\end{prop}
\begin{proof}
	Consider the sign changed invariants $\ogwb^\circ$ with $r = 2$ as in Definition \ref{dfn: sign change}. As a reminder, these are defined by $\ogwb^\circ_{\beta,k} = (-1)^{\frac{k-1}{2}}\ogwb_{\beta,k}$. By Proposition \ref{prop: sign change OWDVV}, the invariants $\ogwb^\circ$ still satisfy the open WDVV equations and the open Gromov-Witten axioms. As these, together with the initial value, were the only inputs for the above induction, the result follows. For invariants with inputs $\Gamma_\diamond$ the result follows from the first statement and the wall-crossing formula \ref{thm:wallcross}.
\end{proof}

\section{The quadric -- even dimension}

\subsection{Relative quantum ring}\label{ssec:relquadev}

In this section, we assume $n>2$ is even. For $n=2$ see Section~\ref{sec:quad2}. We work with the basis
\[
\g_j:= [r^*\oFS]^j,\qquad j=0,\ldots,n,
\]
of $\Hh^*(X,L;\R)$.
We also write $\dl_j:=h^j\in H^*(X;\R)$. Note that in the exact sequence~\eqref{eq:splitLn0} we have $\rho(\g_j)=\dl_j$.
By Lemma~\ref{lm:ip}\ref{it:bez}, the Poincar\'e pairing on $\spa(\{\dl_j\}_{j=0}^n)$ is given by
\[
g_{ij}= \int_X \dl_i\wedge\dl_j
=\delta_{i,n-j}\cdot 2,
\]
where $\delta_{k,l}$ is the Kronecker delta, and the inverse matrix $(g^{ij})_{i,j}:={(g_{ij})_{i,j}}^{-1}$ is given by
$g^{ij}=\delta_{i,n-j}\cdot \frac{1}{2}$.
Thus, the relative quantum product is given by
\[
\mem(\g_i,\g_j)
= \sum_{\substack{0\le m\le n\\ d\ge 0}} \frac{1}{2}T^{d}\GW_d(\dl_i,\dl_j,\dl_m)\g_{n-m},
\qquad i,j\in \{0,\ldots,n\}.
\]

\begin{proof}[Proof of Theorem~\ref{thm:rqhqeven}]
By degree axiom~\eqref{ax_c_deg}, the only nonzero values of $\GW_d(\dl_i,\dl_j,\dl_m)$ occur for $d\le 2$. For $d=2$, the value is only nonzero if $i=j=m=n$. We will not need to compute $\GW_2(\dl_n^{\otimes 3})$ explicitly, although we deduce it below in Remark~\ref{rem:gw2}. For $d=0$, the zero axiom~\eqref{ax_c_zero} says the only nontrivial values are when $i+j+m=n$, in which case $\GW_0(\dl_i,\dl_j,\dl_m)=2$ by Section~\ref{ssec:quadbase}.

For $d=1$, the degree axiom~\eqref{ax_c_deg} says the only nontrivial values occur when $i+j+m=2n$. By Lemma~\ref{lm:gwd1}, the resulting values are
\[
\GW_1(\dl_i,\dl_j,\dl_m)=\begin{cases}
0, & (i,j,m)=(0,n,n),\\
4, & i>0, m=n,\\
8, & i>0, m<n.
\end{cases}
\]
Thus, we have
\[
\mem(\g_i,\g_j)=\begin{cases}
\g_{i+j}, & i+j<n,\\
2T+\g_n, & i+j=n,\\
2T\g_{i+j-n}, & n<i+j<2n.
\end{cases}
\]
To compute the case $i=j=n$, use the associativity of $\mem$. For short, we rite $\g_1^{\,\mem\! k}$ for the $k$-fold product of $\g_1$ under $\mem$. So,
\begin{align}\label{eq:mg2}
\mem(\g_n,\g_n) =& \mem(\g_1^{\,\mem n}-2T, \g_1^{\,\mem n}-2T)\notag \\
=& \g_1^{\,\mem{} 2n}-4T\g_1^{\,\mem n}+4T^2\\
=& \mem(\g_1^{\,\mem (n+1)},\g_1^{\,\mem (n-1)}) - 4T\g_1^{\,\mem n}+4T^2 \notag\\
=& \mem(4T\g_1,\g_1^{\,\mem (n-1)})- 4T\g_1^{\,\mem n}+4T^2 \notag\\
=& 4T^2.\notag
\end{align}
In total,
\[
\mem(\g_i,\g_j)=\begin{cases}
\g_{i+j}, & i+j<n,\\
2T+\g_n, & i+j=n,\\
2T\g_{i+j-n}, & n<i+j<2n,\\
4T^2& i=j=n.
\end{cases}
\]

Similarly to the odd dimensional case, we note that $\{\g_1^{\mem j}\}_{j=0}^n$ is a basis of $QH^*(X,L)$ as a module.
To compute the ring structure we need to spell out the higher powers of $\g_1$. Specifically,
\[
\g_1^{\,\mem (n+1)} = \mem(\g_1,2T+\g_n) = 4T\g_1.
\]
Taking the variables $q,x,$ to correspond to $T,\g_1,$ respectively, we get
the desired result.
\end{proof}

\begin{rem}\label{rem:gw2}
It follows from the computation~\eqref{eq:mg2} that $\GW_2(\dl_n^{\otimes 3})=2[T^2](\mem(\g_n,\g_n))=8$.
\end{rem}

\subsection{Absolute quantum ring}\label{ssec:qhquadev}\label{ssec:abquadev}

In this section, we work with the basis $\{\dl_j\}_{j=0}^n\cup \{\dl_\cir\}$ of $H^*(X;\R)$, where $\dl_j=h^j$ as in Section~\ref{ssec:relquadev} and $\dl_\cir=PD([L])$ as in~\eqref{eq:dL}.
For $i,j\in\{0,\ldots,n,\cir\}$, we have
\[
\dl_i*\dl_j =  \sum_{\substack{0\le m\le n\\ d\ge 0}} \frac{1}{2}T^{d}\GW_d(\dl_i,\dl_j,\dl_m)\dl_{n-m}
+
\sum_{d\ge 0}\frac{(-1)^{n/2}}{2}T^d \GW_d(\dl_i,\dl_j,\dl_{\cir})\dl_{\cir}.
\]

\begin{prop}\label{prop:qhquadev}
The absolute quantum cohomology is given by
\[
QH^*(X)\simeq \altfrac{\R[[q]][x,y]}{I},
\]
with
\[
I=(x^{n+1}-4qx, \quad xy, \quad y^2- (-1)^{n/2+1}4q-(-1)^{n/2}x^n).
\]
\end{prop}
\begin{proof}
By the degree axiom~\eqref{ax_c_deg}, $\GW_d(\dl_i,\dl_j,\dl_{\cir})\ne 0$ implies $d\le 1$. For $d=0$, the zero axiom~\eqref{ax_c_zero} gives
\[
\GW_0(\dl_\cir,\dl_\cir,1)=\int_X\dl_\cir^{\wedge 2} =(-1)^{n/2}2.
\]
Recall the relation $xy=0$ from Proposition~\ref{pr:cohdrlag}. This implies that for $i,j\in\{0,\ldots,n\}$,
$
\GW_0(\dl_\cir,\dl_i,\dl_j) = 0.
$
Finally, by the degree axiom~\eqref{ax_c_deg} $\GW_0(\dl_\cir^{\otimes 3})=0$.

For $d=1$, nontrivial values only appear when $|\dl_i|+|\dl_j|+|\dl_\cir|=2n$. In particular, $(i,j)\ne (\cir, \cir)$.
%Let $\lambda\in A^n(X)$ such that $[\lambda]=\dl_{-1}$.
By equation~\eqref{eq:y2}, % Section~\ref{sssec:dRquadRng}, we have $\lambda^{\wedge 2}=(-1)^{n/2}r^*\oFS^{\wedge n}$. Therefore,
we have
\[
\int_Li^*\dl_\cir =\int_X\dl_\cir^{\wedge 2}  =
(-1)^{n/2}2.
\]
Then $\eta:=\frac{(-1)^{n/2}}{2}\dl_\cir$ satisfies $\int_L\eta=1$.
By the open-closed relation~\eqref{eq:magic}, for any $\xi$ with $\int_L \xi=0$ we have
\[
\ogw_{\beta,1}(\gamma_1,\ldots,\gamma_l).
=
-\GW_\beta(\eta+\xi,\dl_\cir,\gamma_1,\ldots,\gamma_l).
\]
Thus, for $i,j\in\{0,\ldots,n\}$,
\begin{multline*}
\GW_1(\dl_i,\dl_\cir,\dl_j)
= \GW_1(\eta+\dl_{i},\dl_\cir,\dl_j) -\GW_1(\eta,\dl_\cir,\dl_j)=\\
= \ogw_{1,1}(\g_j)-\ogw_{1,1}(\g_j) =0,
\end{multline*}
while by Lemma~\ref{lm:11n},
\[
\GW_1(\dl_\cir,\dl_\cir,\dl_n)
= (-1)^{n/2}2\cdot\GW_1(\eta,\dl_\cir,\dl_n)
= (-1)^{n/2+1}2\cdot\ogw_{1,1}(\g_n)
=(-1)^{n/2+1}4.
\]

In total, we get
\[
\dl_i*\dl_j =
\begin{cases}
\rho(\mem(\g_i,\g_j)),& i,j\in\{0,\ldots,n\},\\
\dl_{\cir}, & i=\cir, j=0,\\
-2T\dl_{\cir},& i=\cir, j=n,\\
(-1)^{n/2+1}2T+(-1)^{n/2}\dl_n,& i=j=\cir,\\
0,& i=\cir, j\ne \cir,0,n.
\end{cases}
\]
Recall that $\dl_n=\dl^{* n}-2T$.
Taking the variables $q,x,y,$ to correspond to $T,\dl_1,\dl_{\cir},$ respectively, we get
the desired result.
\end{proof}

\subsection{Full computation of invariants}

In this section we prove Theorem \ref{thm:evencomputability}. First observe that for $I \subset \{0,1, \dots, n, L\}^{l}$ we have:
\begin{itemize}
	\item $\ogw_{\beta,0}(\Gamma_I) = 0$ by definition, as $H_2(X,L) \cong H_2(X)$.
	\item $\ogw_{\beta,1}(\Gamma_I) = \frac{(-1)^{\frac{3n}{2} + 1}}{2}\GW_d(PD(L),PD(L),\Gamma_I)$ by Theorem \ref{thm:magic}.
	\item $\ogw_{\beta,k}(\Gamma_I) = 0$ for $k\geq2$, also by Theorem \ref{thm:magic}.
\end{itemize}

It will thus suffice to compute all closed Gromov-Witten invariants. We will provide a method to compute these. We first observe:

\begin{lm}
	$\GW_\beta(PD(L)^{\otimes 2s+1},\dl_I) = 0$ for all $\beta$, $s$, and $I \subset \{0,1, \dots, n\}^{l}$.
\end{lm}
\begin{proof}
	Consider the symplectomorphism $\phi: \mathbb{CP}^{n+1} \rightarrow \mathbb{CP}^{n+1}$ given by $\phi(X) = -X$. $\phi$ restricts to a symplectomorphism of $V$, mapping $L$ to itself. $\phi|_{L}$ is the antipodal map. Thus, as $L$ is an even-dimensional sphere, $\phi^*(PD(L)) = -PD(L)$. However $\phi^*(\dl_I) = \dl_I$ as $\phi^*h = h$. As $\phi$ is a symplectomorphism we find \begin{align}
	GW_{\beta}(PD(L)^{\otimes 2s+1},\dl_I)&=GW_{\phi_*\beta}(PD(L)^{\otimes 2s+1},\dl_I)\\
	& = GW_\beta(\phi^*(PD(L)^{\otimes 2s+1}),\phi^*\dl_I) \\
	& = (-1)^{2s+1}GW_\beta(PD(L)^{\otimes 2s+1},\dl_I).
	\end{align}
\end{proof}

The next step is to rule out any invariants with constraints only along $PD(L)$.

\begin{lm}
	\label{lm:at least one not on L}
	$GW_\beta(PD(L)^{\otimes s}) = 0$ for all $\beta$ and $s$.
\end{lm}
\begin{proof}
	For $s$ odd this follows from the previous lemma. So assume $s = 2k$.
	The degree axiom of $\GW$ then implies: \[
	4r + 4r\beta + 4k - 6 = 4kr,
	\]
	Here $n = 2r$. Reducing this equation modulo $4$ shows this is impossible.
\end{proof}

Finally, we can obtain all closed Gromov-Witten invariants.

\begin{lm}
	All Gromov-Witten invariants are determined by the two point invariant $\GW_1(\dl_{n-1},\dl_n) = 4$.
\end{lm}
\begin{proof}
	First observe that by the degree equation there are no non-zero invariants $\GW_\beta(\dl_I)$ with $|I| < 2$, and if $|I| = 2$, the only possibility is $\GW_1(\dl_{n-1},\dl_n)$.
	
	Denote $PD(L)$ by $\dl_L$ and introduce the total order $0 < 1 < \dots < n < L$ on $\{0, 1, \dots, n, L\} := J$. Define the following total order on the set of possible $\beta, I$. Here we always order $I$ to be arranged as $i_1 \leq i_2 \leq \dots \leq i_{|I|}$.
	\begin{itemize}
		\item $(\beta_1, I_1) < (\beta_2, I_2)$ when $\beta_1 < \beta_2$.
		\item $(\beta, I_1) < (\beta,I_2)$ when $|I_1| < |I_2|$.
		\item $(\beta,I_1) < (\beta,I_2)$ when $|I_1| = |I_2|$ and $I_1 < I_2$ in the lexicographic ordering on $J^{|I_1|}$.
	\end{itemize}
	We will now perform induction with respect to this total order. To this end, let $GW_\beta(\dl_I)$ be an invariant we want to compute, with $|I| \geq 3$. Write $I = (i_1, i_2, i_3, \bar{I})$. By the divisor and unit axioms, we may assume $i_1 > 1$ (and hence $i_j >1$ for all $j$). By lemma \ref{lm:at least one not on L}: $i_1 \neq L$.
	We apply the WDVV equation (\ref{thm:WDVV}) with $u = i_{1} - 1$, $v = 1$, $w = i_2$, $y = i_3$ in degree $\beta$, with extra derivatives on $\bar{I}$. This yields:
	\begin{align}
		& \GW_0(\dl_1, \dl_{i_1 - 1}, \dl_{n-i_1})\GW_\beta(\dl_I)g^{i_1,n-i_1} \label{WDVV term 1}\\
		& + \sum_{\mu, \nu \in J} \GW_\beta(\dl_1, \dl_{i_1 - 1}, \dl_{\bar{I}}, \dl_{\mu})\GW_0(\dl_{\nu}, \dl_{i_2},\dl_{i_3})g^{\mu,\nu} \label{WDVV term 2}\\
		& + \sum_{\substack{\beta_1 + \beta_2 = \beta \\ \beta_1, \beta_2 < \beta \\ I_1 \sqcup I_2 = \bar{I} \\ \mu, \nu \in J}} \GW_{\beta_1}(\dl_1, \dl_{i_1 - 1},\dl_\mu, \dl_{I_1})\GW_{\beta_2}(\dl_\nu, \dl_{i_2}, \dl_{i_3},\dl_{I_2}) \label{WDVV term 3}\\
		& = \GW_0(\dl_1, \dl_{i_2}, \dl_{n-i_2-1})\GW_\beta(\dl_{i_2+1},\dl_{i_1 - 1},\dl_{i_3}, \dl_{\bar{I}})g^{i_2+1,n-i_2-1} \label{WDVV term 4}\\
		& + \sum_{\mu, \nu \in J} \GW_\beta(\dl_1, \dl_{i_2}, \dl_{\bar{I}},\dl_\mu)\GW_0(\dl_{\nu}, \dl_{i_1 - 1},\dl_{i_3})g^{\mu,\nu} \label{WDVV term 5}\\
		& + \sum_{\substack{\beta_1 + \beta_2 = \beta \\ \beta_1, \beta_2 < \beta \\ I_1 \sqcup I_2 = \bar{I} \\ \mu, \nu \in J}} \GW_{\beta_1}(\dl_1, \dl_{i_2},\dl_\mu, \dl_{I_1})\GW_{\beta_2}(\dl_\nu, \dl_{i_1 - 1}, \dl_{i_3},\dl_{I_2}) \label{WDVV term 6}\\
	\end{align}
	Now \eqref{WDVV term 1} will determine the invariant we want to compute, as $\GW_0(\dl_1, \dl_{i_1}, \dl_{n-{i_1}}) = 4$ and $g^{i_1,n-i_1} = 1/2$. The first term in \eqref{WDVV term 2} is smaller in the total ordering on indices, as by the divisor axiom we can reduce the number of marked points to $|I| - 1$. All terms in line \eqref{WDVV term 3} are lower down in the ordering, as the degrees are smaller than $\beta$. The second term of \eqref{WDVV term 4} is lower down in the ordering because the the lowest index $i_1$ has been lowered by $1$. Also observe here that we abuse notation by writing $i_2 + 1$ and $n-i_2-1$. If $i_2 = L$, this first term in \eqref{WDVV term 4} vanishes as $\dl_1 \cup \dl_L = 0$. The terms in line \eqref{WDVV term 5} again have fewer interior marked points by the divisor equation, and those in line \eqref{WDVV term 6} have smaller degrees. We have thus determined $\GW_\beta(\dl_I)$ in terms of invariants strictly lower down in the ordering. This finishes the proof by induction.
\end{proof}

\section{The quadric -- dimension 2}\label{sec:quad2}
let $X$ be the quadric hypersurface in $\P^3$ given by $\sum_{j=0}^2z_j^2-z_3^2=0$. Let $L \subset X$ be the real locus. In this section we will compute the relative quantum cohomology of $(X,L)$
\subsection{Topological preliminaries}
Consider the map \begin{align}
	\phi: \mathbb{CP}^1 \times \mathbb{CP}^1 &\rightarrow \mathbb{CP}^3\\
	([x_0:x_1],[y_0:y_1]) &\mapsto [x_0y_0 - x_1y_1: -i(x_0y_0 + x_1y_1): x_1y_0 + x_0y_1: x_0y_1 - x_1y_0]
\end{align}
One then checks that $\phi$ is injective and $im(\phi) = X$. Thus $X \simeq \mathbb{CP}^1 \times \mathbb{CP}^1$. Let $A_1, A_2 \in H_2(X;\mathbb{Z})$ be defined by $A_1 = \phi_*[\mathbb{CP}^1 \times \{[0:1]\}]$ and $A_2 = \phi_*[\{[0:1]\} \times \mathbb{CP}^1]$. Let $l = PD(A_1), l^* = PD(A_2) \in H^2(X;\mathbb{R})$. Then \[
	H^*(X;\mathbb{R}) \simeq \altfrac{\mathbb{R}[l,l^*]}{(l^2, \quad {l^*}^2)}.
\]
Let $L \subset X$ be given by the fixed locus of complex conjugation. Thus $L \simeq S^2$. We now compute $[L] \in H_2(X;\Z)$. An easy computation shows that $|[L] \cap A_i| = 1$, so that $[L] = \pm A_1 \pm A_2$. It remains to compute the signs. In order to do this, we consider $\phi|_L^{-1}: L \rightarrow \mathbb{CP}^1 \times \mathbb{CP}^1$. To compare orientations, it is easiest to precompose this with the (orientation preserving) stereographic projection \begin{align}
\psi: \R^2 &\rightarrow L \subset X\\
(u,v) &\mapsto [2u:2v:u^2+v^2 -1:1+u^2+v^2]
\end{align}
On the image of $\psi$, we have $z_3 \neq z_4$ so that the map $\phi|_{im(\psi)}^{-1}: L \rightarrow \mathbb{CP}^1 \times \mathbb{CP}^1$ given by
\[
[z_0:z_1:z_2:z_3] \mapsto [z_1 + iz_2:z_3 -z_4],[z_3-z_4:iz_2-z_1],
\]
is well defined. The composition with $\psi$ is then given by: \[
(u,v) \mapsto [-(u+iv):1],[1:u-iv].
\]
Now consider the maps $\pi_i: \mathbb{CP}^1 \times \mathbb{CP}^1 \rightarrow \mathbb{CP}^1$. The above then shows that the map $\pi_1 \circ \phi|_{im(\psi)}^{-1} \circ \psi$ is orientation preserving, whereas the map $\pi_1 \circ \phi|_{im(\psi)}^{-1} \circ \psi$ reverses orientation. Thus \[
[L] = A_1 - A_2 \in H_2(X;\Z).
\]
\subsection{Absolute quantum ring}
As shown in~\cite[Example 7.2]{CrauderMiranda} the absolute quantum cohomology of $X$ is given by
\[
QH^*(X;\Lambda_c)\simeq \altfrac{\R[[q_1,q_2]][l,l^*]}{(l^2-q_2, {l^*}^2-q_1)}.
\]
Here
\begin{gather*}
	\L_c=\left\{\sum_{i=0}^\infty a_iT^{\beta_i}\bigg|a_i\in\R,\beta_i\in H_2(X;\Z),\omega(\beta_i)\ge 0,\; \lim_{i\to \infty}\omega(\beta_i)=\infty\right\}.
\end{gather*}
Where we have taken the identification \[
\L_c \simeq \R[[q_1,q_2]]
\]
via the map $T^{A_1} \mapsto q_1$, $T^{A_2} \mapsto q_2$. However, to compare $QH^*(X)$ with $QH^*(X,L)$ we need to work over a Novikov ring with exponents in $H_2(X,L;\Z)$. As $[L] = A_1 - A_2$, and $H_1(L;\Z) = 0$ this amounts to identifying $q_1 = q_2$. We additionally change basis to $u = l + l^*$, $v = l - l^*$ to obtain:
	\[
	QH^*(X)\simeq \altfrac{\R[[q]][u,v]}{(u^2+v^2-4q, \quad uv)}.
	\]
Note that this latter result agrees with the quantum cohomology for the general even dimensional quadric.

\subsection{Relative quantum ring}
As $[L] \neq 0$, we have $QH^*(X;L) \subset QH^*(X)$ given by the kernel of integration over $L$, $\i_L$. As a $\R[[q]]$-module, $QH^*(X,L)$ is thus spanned by $1, u$ and $l\cup\la$. From the previous section, we have that $u \star u = 2q + 2 l\cup \la$. We thus find:
\[
QH^*(X,L)\simeq \altfrac{\R[[q]][u]}{(u^3 - 4uq)}.
\]
We have thus proven Theorem~\ref{thm:rqhq2}.

\subsection{Full computation of invariants}
$(1,0)=[\ml]$, $(0,1)=[\mla]$.

\begin{lm}\label{lm:gw2nn-1}
	The basic Gromov-Witten invariants are given by
	\[
	\GW_{(d_1A_1 + d_2A_2)}(l\cup\la, \a) =
	\begin{cases}
		1, & (d_1,d_2, \a) \in \{(1,0,l), (0,1,\la)\},\\
		1, & (d_1,d_2,\a) = (0,0,l\cup\la),\\
		0, & \text{otherwise}.
	\end{cases}
	\]
\end{lm}
\begin{proof}[Proof of Theorem~\ref{thm:q2computability}]
	As $H^2(X)$ generates $H^*(X)$ as an algebra, the result for closed invariants follows from \cite{KontsevichManin}. As $[L] \neq 0$, the result for open Gromov-Witten invariants follows from the result for closed invariants and Theorem \ref{thm:magic}.
\end{proof}

\appendix
\section{Initial value}\label{sec:app}
Our goal here is to prove Lemma~\ref{lm:init}.

\subsection{Setting up the computation}
For any space, denote by $pt$ the map from it to a point.
We use the notation of~\cite{ST2,ST3} to write $\ogw_{1,3}$ as
\begin{align}
(\d_s^3\q^{\beta_1;\gamma,b}_{-1,0})|_{s=t_j=0}
= &
(\d_s^3\q^{\beta_1;0,b}_{-1,0})|_{s=t_j=0} \notag\\
=&
\langle\q^{\beta_1}_{1,0}(\d_s b,\d_s b),\d_s b\rangle \notag\\
=&
(-1)^{n+1} pt_*((evb_0)_*(evb_1^*\d_s b\wedge evb_2^*\d_s b)\wedge \d_s b) \notag \\
=&
(-1)^{1+n\cdot\rdim evb_0} pt_*(\d_s b \wedge (evb_0)_*(evb_1^*\d_s b\wedge evb_2^*\d_s b)) \notag \\
=&
pt_*(evb_0^*\d_s b\wedge evb_1^*\d_s b\wedge evb_2^*\d_s b), \label{eq:ds3}
\end{align}
where $b$ is (the) point-like bounding chain guaranteed by Theorem~2 of~\cite{ST2}.

Let $Y = L \times L \times L \setminus \triangle$ and let $Z = (evb_0 \times evb_1 \times evb_2)^{-1}(Y) \subset \M_{3,0}(\beta_1).$
Observe that
\[
g = (evb_0 \times evb_1 \times evb_2)|_Z : Z \to Y
\]
is a covering map. Indeed, the fibre of $g$ over a point $(x_1,x_2,x_3) \in Y$ can be identified with the pair of real oriented degree-2 curves in $L$ passing through $x_1$, $x_2$, and $x_3.$
We claim that for the unique choice of (relative) spin structure with $w_\s = 0$, the degree of $g$ is $2.$ Indeed, the two points in the preimage of a point in $Y$ are conjugate disks of degree $1$ with three marked points. Proposition~5.1 of~\cite{SolomonThesis} shows that the covering transformation that interchanges these two disks is orientation preserving:
\begin{align*}
s^+(1,3,0)  &\equiv  \frac{2n(2n+1)}{2}+n+n+3+0 + (n+1)\sum_{a}\frac{(k_a-1)(k_a-2)}{2}\\
&\equiv n(2n+1) + 1
\equiv 0 \pmod 2.
\end{align*}
So, the degree is $\pm 2.$
We postpone the proof of the following proposition until later in this section.
\begin{prop}\label{prop:gplus}
The unique spin structure on $L$ is such that $g$ has degree $(-1)^{\frac{(n-2)(n-3)}{2}}2$.
\end{prop}
%Given this result, we deduce that for the above mentioned choice of relative spin structure, the degree of $g$ is $(-1)^{\frac{(n-2)(n-3)}{2}}2$.

For $i = 1,2,3,$ let $p_i: L \times L \times L \to L$ be the projections. Since $Y \subset L \times L \times L$ and $Z \subset \M_{3,0}(\beta_1)$ are open dense subsets, using that $b$ is point-like, we obtain
\begin{align*}
pt_*((evb_0)^*\d_s b\wedge (evb_1)^*\d_s b\wedge(evb_2)^*\d_s b) & = \int_Zg^*((p_1^* \d_s b \wedge p_2^* \d_s b \wedge p_3^*\d_s b)|_Y) \\
&= (-1)^{\frac{(n-2)(n-3)}{2}}2 \int_Y p_1^* \d_sb \wedge p_2^* \d_s b \wedge p_3^* \d_s b\\
&= (-1)^{\frac{(n-2)(n-3)}{2}}2 \int_{L\times L \times L} p_1^* \d_sb \wedge p_2^* \d_s b \wedge p_3^* \d_s b\\
&=(-1)^{\frac{(n-2)(n-3)}{2}}2\left (\int_L \d_s b \right)^3 \\
&= (-1)^{\frac{(n-2)(n-3)}{2}}2.
\end{align*}
Combining this calculation with equation~\eqref{eq:ds3} we get that the value of $\ogwb_{1,3}$ is $(-1)^\frac{(n-2)(n-3)}{2}2$.
%However, a different choice would alter the sign of $\int_L\d_s b$, so the change of sign in $\deg g$ will cancel out with the change of sign in $\big(\int_L \d_s b\big)^3$.
%Thus, the sign of $\ogwb_{1,3}$ in fact does not depend on the choice of the relative spin structure.

\begin{rem}
The choice of just the orientation on $L$ has no effect on the value, since a change of orientation would contribute $(-1)^4=+1$ to the sign of $\left(\int_L\d_sb\right)^3$.
\end{rem}

It is left for us to prove Proposition~\ref{prop:gplus}. At a point
\[
\uu=[u,\vec{z}=(1,i,-1),\vec{w}=()]\in \M_{3,0}(\beta_1),
\]
we need to compute the sign of the evaluation map
\[
dg_{\uu}:T_{\uu}\M_{3,0}(\beta_1)=\g(u^*TX,u^*TL)
\lrarr T_{u(1)}L\oplus T_{u(i)}T\oplus T_{u(-1)}L.
\]
In section \ref{sec: sign is independent of map} we prove:
\begin{lm}
	\label{lm: sign is independent of map}
	The sign of $dg_{\uu}$ is independent of the choice of $\uu$.
\end{lm}
Thus, we can pick a convenient map $\uu_0$ and compute the sign $dg_{\uu_0}$. We outline this computation. First, consider the Cauchy-Riemann boundary value problem \[
\E=(E,F)
\]
given by the trivial line bundle over $(D,\d D)$ with boundary conditions
\[
F_z=z\R,\qquad z\in\d D.
\]
The Cauchy-Riemann boundary problem $\E$ has a canonical orientation and spin structure, $\mathfrak{p}_\E$, given by trivialising the bundle $F_z$ by division by $z$. This then equips $\g(\E)$ with a natural orientation (see e.g. \cite{SolomonThesis}). In section \ref{sec: sign of evaluation} we show:
\begin{lm}
	\label{lm: sign of evaluation}
	The canonical orientation and spin structure on $\E$ equip $\g(\E)$ with an orientation such that \[
	\Gamma(\E) \xrightarrow{Ev} F_1 \oplus F_i \oplus F_{-1}
	\]
	is orientation preserving.
\end{lm}
We then show that the problem we are interested in decomposes into copies of these line bundles: \begin{lm}
	\label{lm: decomposition into line bundles}
	There exists an isomorphism of Cauchy-Riemann boundary value problems \[
	(u_0^*TX,u^*TL) \xrightarrow{\varphi} \E^{\oplus n}.
	\]
	Moreover, the map $\varphi$ does not respect spin structures.
\end{lm}
The proof of the above lemma is given in section \ref{sec: decomposition into line bundles}. Then consider the commutative diagram:
\[
\xymatrix{
\g(u^*TX,u^*TL)\ar[rr]^{dg_{\uu}}\ar[d]^{\varphi_*} & & T_{u(1)}L\oplus T_{u(i)}L\oplus T_{u(-1)}L\\
\g(\E^{\oplus n})\ar[d] & & \\
\g(\E)^{\oplus n}\ar[r]^(0.35){Ev} & (F_1\oplus F_{i}\oplus F_{-1})^{\oplus n}\ar[r]^{\sigma} & F_1^{\oplus n}\oplus F_i^{\oplus n}\oplus F_{-1}^{\oplus n}\ar[uu]_{\tau}
}
\]
Here $\tau$ is the map induced by restricting $\varphi$ to the fibres over $1, i$ and $-1$. To recover the sign of the arrow on the top row, we compute the signs of all the intermediate steps. From Lemma \ref{lm: decomposition into line bundles} we have that $sgn(\varphi_*) = 1 + sgn(\varphi|_{TL})$, where the $1$ comes from the fact that it doesn't respect the spin structure, and $sgn(\varphi|_{TL})$ is the sign coming from whether or not the map $\varphi: u^*TL \rightarrow F^{\oplus n}$ respects orientations. Next, note that the natural map $\g(\E^{\oplus n}) \rightarrow \g(\E)^{\oplus n}$ is orientation preserving by \cite[Lemma~2.16]{SolomonThesis}. Then, we have the map $Ev$ which, as shown in Lemma \ref{lm: sign of evaluation}, is orientation preserving. The map $\sigma$ is just a permutation, with sign $sgn(\sigma)=\frac{n(n-1)}{2}.$ Finally, we have the map $\tau$, which is induced by $\varphi$. It has sign $sgn(\varphi|_{TL})^3 = sgn(\varphi|_{TL})$. The overall sign is thus \[
1 + \frac{n(n-1)}{2} \equiv \frac{(n-2)(n-3)}{2},
\]
as required. Note in particular that we don't need to compute $sgn(\varphi|_{TL})$.

\subsection{Proof of Lemma \ref{lm: sign is independent of map}}
\label{sec: sign is independent of map}
Showing $sgn(dg_{\uu})$ is independent of the choice of $\uu$.

Let $Q : \C^{n+2} \to \C$ be the quadratic form given by
\[
Q(z) = \sum_{j=0}^n z_j^2-z_{n+1}^2.
\]
Let
\[
\Xh:=\{Q(z)=0\}\subset \C^{n+2}
\]
and the real locus
\[
\Lh\subset \R^{n+2}.
\]
Let $e_0,\ldots,e_{n+1} \in \C^{n+2}$ denote the standard basis and let $e_0^*,\ldots,e_{n+1}^*$ denote the dual basis.
Let
\[
\langle \cdot, \cdot \rangle_{n+1,1} =  \sum_{j = 0}^n e_j^* \otimes e_j^* - e_{n+1}^* \otimes e_{n+1}^*
\]
denote the standard bilinear form of signature $(n+1,1).$
Let $P \subset \R^{n+2}$ be an oriented $3$-dimensional linear subspace such that $\langle \cdot, \cdot \rangle_{n+1,1}|_P$ has signature $(2,1).$ Let $P^\C  = P \otimes \C \subset \C^{n+2}.$ Observe that $\bP(P \cap \Lh)$ is isomorphic as an algebraic variety to $\RP^1.$ On the other hand, the projection map $P_L := \Lh \cap P \cap \{ z_{n+1} = 1\} \to \bP(P \cap \Lh)$ is a diffeomorphism. We equip $\bP(P \cap \Lh)$ with the induced orientation, where $P_L$ is equipped with the following orientation. A tangent vector $v$ to $P_L$ is positively oriented if there is a vector $w \in P$ such that the directional derivative $\d_wQ$ is positive, and vector $x \in P$ such that $\langle x,e_{n+1}\rangle_{n+1,1}$ is positive and $w,v,x$ is an oriented basis for $P.$ Equip $\RP^1$ with the orientation induced by viewing it as the boundary of the hemisphere in $\P^1$ containing $[1:\i]$. Let $\tilde u_P^\R : \RP^1 \to \bP(P \cap \Lh)$ be an orientation preserving algebraic isomorphism. So, $\tilde u_P^\R$ is uniquely determined by $P$ up to reparameterization by $PSL_2(\R).$ Let $\tilde u_P : \P^1 \to \bP(P^\C \cap \Xh)$ be the complexification of $\tilde u_P^\R.$
Let
\[
\psi : D \to \P^1
\]
be given by
\[
\psi(z) = [\i z + \i : z-1], \qquad z \in \C, \; |z| \leq 1.
\]
Let $u_P : (D,\d D) \to (X,L)$ be the holomorphic map given by $u_P = \tilde u_P \circ \psi.$

Let $Gr(2,n+2)^{(2,1)}$ denote the Grassmanian of oriented subspaces $P \subset \R^{n+2}$ as above. Let $SO(n+1,1)_+$ denote the group of linear automorphisms $l : \R^{n+2} \to \R^{n+2}$ preserving the bilinear form $\langle \cdot, \cdot \rangle_{n+1,1}$ and satisfying
\[
\det l = 1, \qquad \langle  l(e_{n+1}), e_{n+1} \rangle_{n+1,1} > 0.
\]
The group $SO(n+1,1)_+$ is connected. It acts on $Gr(2,n+2)^{(2,1)}$ and this action is transitive for $n > 1.$ Moreover, $SO(n+1,1)_+$ acts on $(X,L)$ by holomorphic symplectomorphisms and thus also acts on the moduli spaces $\M_{k,l}(\beta).$
\begin{lm}\label{lm:Gr}
The map $Gr(2,n+2)^{(2,1)} \to \M_{3,0}(\beta_1)$ given by $P \mapsto u_P$ is an $SO(n+1,1)_+$ equivariant diffeomorphism.
\end{lm}
Then, as the map $\M_{3,0}(\beta_1) \rightarrow L^3$ is also $SO(n+1,1)_+$ equivariant, and $SO(n+1,1)_+$ is connected (and hence acts in an orientation preserving way), we find that the sign of the evaluation map is independent of the choice of map.

\subsection{Proof of Lemma~\ref{lm: decomposition into line bundles}}
\label{sec: decomposition into line bundles}
Decomposing $(u^*TX,u^*TL)$ into line bundles.

Let
\[
P_0=\spa\{e_0,e_1,e_{n+1}\}.
\]
We compute $u_{P_0}$ explicitly. Throughout, we shorten the subscript $P_0$ to $0.$ $P_0$ has associated map $u_0: (D,\d D)\lrarr (X,L)$ given by:
\[
u_0(z) = [z^2+1: \i (1-z^2): 0 : \cdots : 0 : 2z].
\]
After normalisation, the restriction to $\d D$ is given by: \[
u_0(\theta) =  [\cos \theta: \sin \theta: 0 :\cdots : 1]
\]
As the map $u_0$ is non-singular, the derivative $Du(\frac{\d }{\d z})$ defines an nowhere vanishing holomorphic section of $u_0^*TX$, call this $\xi_1$. We will now extend this to a basis of sections of $u_0^*TX$. To this end, for $j = 2, \dots, n$, let $\xi_j$ be the vector field given by $(0,0,\dots,0,\frac{1}{z},0,\dots, 0) = \frac{1}{z}\partial_{x_j}$ in the affine chart $z_{n+1} \neq 0$, with coordinates $x_0, \dots, x_n$. We then have:
\begin{lm}
	The vector fields $\xi_j$ extend to holomorphic sections $D \rightarrow u_0^*TX$. Moreover, they are nowhere vanishing.
\end{lm}
\begin{proof}
	The unique extension can be computed by instead going to the affine patch $z_0 \neq 0$. Using the coordinates $y_1, \dots, y_{n+1}$ in this patch, we find that $\xi_j = \frac{2}{1+z^2}\partial_{y_j}$. Thus, $\xi_j$ extends uniquely over the apparent pole at $z= 0$, and is indeed nowhere vanishing.
\end{proof}
Let $(E_j,F_j)$ be the line sub-bundle of $(u_0^*TX, u_0^*TL)$ generated by $\xi_j$. That is, we set $E_j = span_\C \xi_j$, and the real boundary condition is given by $F_j := E_j \cap u_0^*TL$.
\begin{lm}\label{lm:frameconstr}\leavevmode
	The map \[
	\varphi: \E^{\oplus n} \rightarrow (u_0^*TX,u_0^*TL)
	\]
	given by \[
	\varphi_z(\lambda_1, \dots, \lambda_n) = -i\lambda_1 \xi_1 +\lambda_2 \xi_2 + \dots + \lambda_n \xi_n
	\]
	is an isomorphism.
\end{lm}
\begin{proof}
	Note first that if the map is well-defined, it is an isomorphism, as the sections $\xi_j$ are linearly independent (and thus form a basis). It thus remains to be shown that the map respects the real boundary conditions, and is thus well-defined.
	For the first component, we do the computation in the affine patch $\{z_{n+1} = 1\}$. The map $u_0: D \setminus \{ 0 \} \rightarrow Q \subset \mathbb{C}^{n+1}$ is given by \[
	u_0(z) = (\frac{z^2 + 1}{2z}, i\frac{1-z^2}{2z}, 0, \dots, 0).
	\]
	One can check that indeed, if $|z| = 1$, the quantities $\frac{z^2 + 1}{2z}$, and $i\frac{1-z^2}{2z}$ are both real numbers. In fact, if $z = e^{i\theta}$ they are $\cos \theta$ and $\sin \theta$ respectively. The derivative can be computed as \[
	Du(\frac{\d }{\d z}) = (\frac{z^2 - 1}{2z^2}, -i\frac{1+z^2}{2z^2}, 0, \dots, 0)
	\]
	Note that this lies indeed in the tangent bundle of $Q$, which is given by those $(p,v) \in \C^{n+1} \oplus \C^{n+1}$ with $p \cdot v = 0$, where $\cdot$ denotes the complex linear extension of the dot product on $\R^{n+1}$.
	Then, to show that $\varphi_1$ respects the real boundary conditions, note that on the left, it is given by $F_z = z \R$. On the right it is given by $span_\C \xi_1 \cap u_0^*TL$. Thus, let $\lambda_1 = \alpha_1 z$ for $\alpha_1 \in \R$. We compute \[
	-i \alpha_1 z \xi_1 = \alpha_1 (i\frac{1 - z^2}{2z}, \frac{1+z^2}{2z}, 0, \dots, 0),
	\]
	and, indeed, as noted before, these two quantities are real numbers, so that $-i \alpha_1 z \xi_1 \in u_0^*TL$ as required. Next, for $j \geq 2$ we have \[
	\xi_j = (0, 0, \dots, \frac{1}{z}, 0, \dots, 0)
	\]
	Observe that if $\lambda_j = \alpha_j z$ for $\alpha_j \in \R$, we have $\alpha_j z \xi_j = (0, \dots, 0, \alpha_j, 0, \dots, 0) \in u_0^*TL$ as required.
\end{proof}

We are now ready to compute the effect on spin structures. The spin structure on $u_0^*TL$ is the unique spin structure obtained as follows: take the loop $u_0|_{\d D} \subset L$. As $L$ is simply connected, one can choose a disk $D$ in $L$ bounding this loop. The bundle $TL|_D$ is trivial (unique up to homotopy) as $D$ is contractible. Restricting this trivialisation to $\d D$ yields a trivialisation of $u_0|_{\d D}^* TL$. Note that as $w_2(L) = 0$, so any two choices of $D$ yield homotopic trivialisations of $u_0|_{\d D}^* TL$. Now equip the resulting trivial bundle $u_0^*TL$ bundle with the trivial spin structure (i.e. the one that in this trivialisation is just the trivial $spin(n))$ bundle over $S^1$. We will explicitly write down such a trivialisation.
\begin{lm}
	The natural trivialisation of $u_0^*TL$ coming from a disk bounding $u_0|_{\d D}$ has a framing homotopic to the one given by $\{\psi_j\}_{j=1, \dots, n}$ with \begin{align*}
	\psi_1 &= (\sin^2 \theta, -\cos \theta \sin \theta, \cos \theta, 0 , \dots, 0)\\
	\psi_2 &= (-\cos \theta \sin \theta, \cos^2 \theta, \sin \theta, 0, \dots, 0)\\
	\psi_j &= (0, 0, 0, \dots, 0, 1, 0, \dots, 0) \text{ for } j \geq 3
	\end{align*}
	Here $\psi_j$ has a non-zero entry only in entry $j+1$. The vector fields are expressed by viewing $L \subset \R^{n+1}$ as the unit sphere.
\end{lm}
\begin{proof}
	Using the same notation as in the rest of the section, the coordinates on $\R^{n+1}$ are $x_0, \dots, x_n$. Let the disk $\mathbb{D}$ be given by $\{x_0^2 + x_1^2 + x_2^2 = 1, x_2 \geq 0| x_i \in \R\}$. Identify $\mathbb{D}$ with $D$ using stereographic projection. Note that we have $u_0|_{S^1}^*TL \simeq T\mathbb{D}|_{\d \mathbb{D}} \oplus \R^{n-2}$, where the $\R^{n-2}$ factor is given by the vector fields $\psi_j$ for $j \geq 3$. Now, trivialise $TD$ with the standard unit vector fields $\partial_x$, $\partial_y$. Pulling these back under stereographic projection and restricting to the boundary yields the vector fields $\psi_1$ and $\psi_2$.
\end{proof}

Next, we look at the trivialisation induced by the isomorphism $(u_0^*TX, u_0^*TL) \simeq \E^{\oplus n}$. The boundary condition for $\E$ is given by $F_z = z\R$, and is trivialised by division by $z$. The spin structure is again given by taking the trivial spin structure with respect to this trivialisation. We have:
\begin{lm}
	The trivialisation of $u_0^*TL$ coming from the isomorphism with $\E^{\oplus n}$ has a framing homotopic to the one given by $\{ \chi_j\}_{j=1, \dots, n}$ with \begin{align*}
		\chi_1 &= (\sin \theta, -\cos \theta , 0, 0 , \dots, 0)\\
		\chi_j &= (0, \dots, 0, 1, 0, \dots, 0) \text{ for } j \geq 2
	\end{align*}
	where $\chi_j$ has a non-zero entry only in position $j+1$.
\end{lm}
\begin{proof}
	The sections $\chi_j$ are induced by the sections $\xi_j$ from before. Concretely, recall that \[
	\xi_1 = (\frac{z^2 - 1}{z^2}, -i\frac{1+z^2}{z^2}, 0, \dots, 0).
	\]
	The morphism $\varphi_1: F_z \rightarrow u_0^*TL$ is given by $\alpha z \mapsto -i\alpha z \xi_1$, which after substituting $z = e^{i\theta}$ indeed yields $\chi_1$. For $j\geq 2$ the morphism $\varphi_j: F_z \rightarrow u_0^*TL$ is given by $\alpha z \mapsto \alpha \xi_j = \alpha (0, \dots,0, 1, 0, \dots, 0)$
thus the induced section is indeed $\chi_j$.
\end{proof}
Finally, we have:
\begin{lm}
	The trivialisations induced by the basis $\psi_j$ is not homotopic to the one induced by $\chi_j$.
\end{lm}
\begin{proof}
	The change of basis matrix is given by \begin{equation*}
		\begin{pmatrix}
			\sin \theta & \cos \theta & 0\\
			-\cos \theta & \sin \theta & 0\\
			0 & 0 & I_{n-2}
		\end{pmatrix}
	\end{equation*}
	which indeed defines a non-trivial loop in $SO(n)$.
	\end{proof}
The conclusion is that the induced spin structures are opposed.

\subsection{Proof of Proposition~\ref{prop:gplus}}
\label{sec: sign of evaluation}
Computing the sign of $Ev$.

It is known that global holomorphic sections of $\E$ over $(D,\d D)$ are described as follows:
\[
\g(\E) = \{a_0 +a_1z+\bar{a}_0z^2 \;|\; a_1\in \R\}.
\]
Each copy of $\E\subset T_\uu\M_{3,0}(1)$ gives rise to maps
\begin{gather*}
evb^1: \g(\E)\lrarr F_1\simeq \R,\\
evi^0: \g(\E)\lrarr E_0\simeq \C,
\end{gather*}
given by evaluation at $1\in \d D$ and $0\in int(D)$, respectively. Lemma A.2.1 of~\cite{Smith} says that in fact this determines an orientation-preserving isomorphism
\[
\theta= evb^1\oplus evi^0: \g(\E) \stackrel{\sim}{\lrarr} \R\oplus \C.
\]
Using $\theta$, we find that the following sections form a positively oriented basis for $\g(\E)$, when taken in this order:
\begin{gather*}
\xi^1(z) = 1-2z+z^2,\\
\xi^2(z) = i-iz^2,\\
\xi^3(z) = z.
\end{gather*}
Evaluating $\xi^j$ at $1, i,$ and $-1,$ we see that
\[
\begin{pmatrix}
\xi^1(1) & \xi^2(1) & \xi^3 (1)\\
\xi^1(i) & \xi^2(i) & \xi^3(i)\\
\xi^1(-1) & \xi^2(-1) & \xi^3(-1)
\end{pmatrix}
=
\begin{pmatrix}
0 & 0 & 1\\
-2i & 2i & i\\
4 & 0 & -1
\end{pmatrix}.
\]
By definition, the fibres of the boundary condition $F$ of $\E$ at $z=1,i,-1,$ are, respectively
\[
F_1 = \R, \quad
F_i = i\R,\quad
F_{-1} = -\R.
\]
Thus,
\[
evb:\g(\E)\stackrel{\sim}{\lrarr} \R\oplus i\R\oplus -\R
\]
is represented by the matrix
\[
M=
\begin{pmatrix}
0& 0& 1\\
-2& 2& 1\\
-4& 0& 1
\end{pmatrix},
\]
which has determinant
\[
\det M =
1\cdot \begin{vmatrix}
-2& 2\\
-4&0
\end{vmatrix}=
8>0.
\]
So, $Ev$ is represented a block-diagonal matrix where the diagonal consists of $n$ copies of the block $M$ and satisfies
\[
sgn(Ev)=sgn(\det M)^n=1.
\]

\bibliography{../../bibliography_exp}

\providecommand{\bysame}{\leavevmode\hbox to3em{\hrulefill}\thinspace}
\providecommand{\MR}{\relax\ifhmode\unskip\space\fi MR }
\providecommand{\MRhref}[2]{%
  \href{http://www.ams.org/mathscinet-getitem?mr=#1}{#2}
}
\providecommand{\href}[2]{#2}
\begin{thebibliography}{10}

\bibitem{Beauville}
A.~Beauville, \emph{Quantum cohomology of complete intersections}, Mat. Fiz.
  Anal. Geom. \textbf{2} (1995), no.~3-4, 384--398.

\bibitem{CieliebakGoldstein}
K.~Cieliebak and E.~Goldstein, \emph{A note on the mean curvature, {M}aslov
  class and symplectic area of {L}agrangian immersions}, J. Symplectic Geom.
  \textbf{2} (2004), no.~2, 261--266.

\bibitem{CrauderMiranda}
B.~Crauder and R.~Miranda, \emph{Quantum cohomology of rational surfaces}, The
  moduli space of curves ({T}exel {I}sland, 1994), Progr. Math., vol. 129,
  Birkh\"{a}user Boston, Boston, MA, 1995, pp.~33--80, \href
  {http://dx.doi.org/10.1007/978-1-4612-4264-2\_3}
  {\path{doi:10.1007/978-1-4612-4264-2\_3}}.

\bibitem{FultonPandharipande}
W.~Fulton and R.~Pandharipande, \emph{Notes on stable maps and quantum
  cohomology}, Algebraic geometry---{S}anta {C}ruz 1995, Proc. Sympos. Pure
  Math., vol.~62, Amer. Math. Soc., Providence, RI, 1997, pp.~45--96, \href
  {http://dx.doi.org/10.1090/pspum/062.2/1492534}
  {\path{doi:10.1090/pspum/062.2/1492534}}.

\bibitem{KontsevichManin}
M.~Kontsevich and Y.~Manin, \emph{Gromov-{W}itten classes, quantum cohomology,
  and enumerative geometry}, Comm. Math. Phys. \textbf{164} (1994), no.~3,
  525--562.

\bibitem{LercheVafaWarner}
W.~Lerche, C.~Vafa, and N.~P. Warner, \emph{Chiral rings in {$N=2$}
  superconformal theories}, Nuclear Phys. B \textbf{324} (1989), no.~2,
  427--474, \href {http://dx.doi.org/10.1016/0550-3213(89)90474-4}
  {\path{doi:10.1016/0550-3213(89)90474-4}}.

\bibitem{MS}
D.~McDuff and D.~Salamon, \emph{{$J$}-holomorphic curves and symplectic
  topology}, second ed., American Mathematical Society Colloquium Publications,
  vol.~52, American Mathematical Society, Providence, RI, 2012.

\bibitem{RuanTian0}
Y.~Ruan and G.~Tian, \emph{A mathematical theory of quantum cohomology}, Math.
  Res. Lett. \textbf{1} (1994), no.~2, 269--278.

\bibitem{RuanTian}
Y.~Ruan and G.~Tian, \emph{A mathematical theory of quantum cohomology}, J.
  Differential Geom. \textbf{42} (1995), no.~2, 259--367.

\bibitem{Smith}
J.~Smith, \emph{Discrete and continuous symmetries in monotone {F}loer theory},
  Selecta Math. (N.S.) \textbf{26} (2020), no.~3, Paper No. 47, 65, \href
  {http://dx.doi.org/10.1007/s00029-020-00575-5}
  {\path{doi:10.1007/s00029-020-00575-5}}.

\bibitem{SolomonThesis}
J.~P. {Solomon}, \emph{{Intersection theory on the moduli space of holomorphic
  curves with Lagrangian boundary conditions}}, arXiv e-prints (2006), \href
  {http://arxiv.org/abs/math/0606429} {\path{arXiv:math/0606429}}.

\bibitem{ST3}
J.~P. {Solomon} and S.~B. {Tukachinsky}, \emph{{Relative quantum cohomology}},
  to appear in J. Eur. Math. Soc., \href {http://dx.doi.org/10.4171/JEMS/1337}
  {\path{doi:10.4171/JEMS/1337}}.

\bibitem{ST2}
J.~P. Solomon and S.~B. Tukachinsky, \emph{Point-like bounding chains in open
  {G}romov-{W}itten theory}, Geom. Funct. Anal. \textbf{31} (2021), no.~5,
  1245--1320, \href {http://dx.doi.org/10.1007/s00039-021-00583-3}
  {\path{doi:10.1007/s00039-021-00583-3}}.

\bibitem{Witten1}
E.~Witten, \emph{Topological sigma models}, Comm. Math. Phys. \textbf{118}
  (1988), no.~3, 411--449.

\end{thebibliography}
\bibliographystyle{../../amsabbrvcnobysame}

\end{document}